\RequirePackage{scrlfile}
\newlength{\bibsep}
\documentclass[3p, nonatbib, final]{elsarticle}

\usepackage[utf8]{inputenc}
\usepackage[british]{babel}

\usepackage{csquotes}
\makeatletter
\let\c@author\relax
\let\bibfont\relax
\let\bibsep\relax
\makeatother
\usepackage[
    sorting=none,
    style=numeric-comp,
    doi=true,
    url=true,
    maxbibnames=5,
    backend=biber,
    giveninits=true,
    isbn=false,
    ]{biblatex}
\bibliography{literature}
\renewcommand*{\bibfont}{\normalfont\footnotesize}
\usepackage{graphicx}
\usepackage[export]{adjustbox}

\usepackage{amssymb,amsmath, amsthm}
\usepackage{mathtools}
\usepackage{commath}
\usepackage{nicefrac}
\usepackage{array}
\usepackage{multirow}
\newcommand\bmmax{2}
\usepackage{bm}
\usepackage{calrsfs}
\usepackage{stmaryrd}
\DeclareMathAlphabet{\pazocal}{OMS}{zplm}{m}{n}
\DeclareMathAlphabet{\pazocalbf}{OMS}{cmsy}{b}{n}
\usepackage{etoolbox}
\preto\subequations{\ifhmode\unskip\fi}

\usepackage[dvipsnames]{xcolor}
\usepackage[colorlinks=true]{hyperref}
\usepackage[nameinlink]{cleveref}

\crefname{algocf}{algorithm}{algorithms}

\usepackage{enumitem}
  \setlist[enumerate]{nosep, topsep=0pt, wide = 1em, leftmargin=*}
  \setlist[itemize]{nosep, topsep=3pt, wide = 1em, leftmargin=*}

\usepackage[linesnumbered, ruled]{algorithm2e}

\newtheorem{form}{Formulation}

\theoremstyle{remark}
\newtheorem{remark}{Remark}

\def\p{\partial}
\def\eps{\varepsilon}
\let\phivar\phi
\def\phi{\varphi}

\newcommand{\bs}[1]{{#1}}
\newcommand{\bg}{\bs{g}}
\newcommand{\hcalA}{\hat{\cal A}}

\newcommand{\RR}{\mathbb{R}}
\newcommand{\PP}{\mathbb{P}}

\newcommand{\xb}{\bm{x}}
\newcommand{\xbref}{\hat{\xb}}
\newcommand{\nb}{\bm{n}}
\newcommand{\db}{\bm{d}}

\newcommand{\CR}{\pazocal{C}}
\newcommand{\CRh}{\CR_h}
\renewcommand{\O}{\Omega}
\newcommand{\Oref}{\hat{\O}}

\newcommand{\FL}{\pazocal{F}}
\newcommand{\FLref}{\hat{\pazocal{F}}}
\newcommand{\SO}{\pazocal{S}}
\newcommand{\SOref}{\hat{\pazocal{S}}}
\newcommand{\IN}{\pazocal{I}}
\newcommand{\INref}{\hat{\pazocal{I}}}

\newcommand{\phiref}{\hat{\bm{\phivar}}}
\newcommand{\xiref}{\hat{\xi}}
\newcommand{\psiref}{\hat{\bm{\psi}}}

\newcommand{\fb}{\bm{f}}
\newcommand{\fbref}{\hat{\fb}}
\newcommand{\gb}{\bm{g}}
\newcommand{\gbref}{\hat{\gb}}

\newcommand{\ub}{\bm{u}}
\newcommand{\ubref}{\hat{\ub}}
\newcommand{\ubh}{\ub_h}
\newcommand{\vb}{\bm{v}}
\newcommand{\vbref}{\hat{\vb}}
\newcommand{\vbh}{\vb_h}
\newcommand{\pref}{\hat{p}}
\newcommand{\wb}{\bm{w}}

\newcommand{\sigb}{\bm{\sigma}}
\newcommand{\sigbref}{\hat{\sigb}}
\newcommand{\eb}{\bm{e}}

\newcommand{\Vb}{\bm{V}}
\newcommand{\Vbref}{\hat{\Vb}}
\newcommand{\Vbh}{\Vb_h}
\newcommand{\Qh}{Q_h}

\newcommand{\dn}{\partial_{\nb}}
\newcommand{\dnF}{\partial_{\nb_F}}
\newcommand{\jump}[1]{\llbracket#1\rrbracket}

\newcommand{\Th}{\mathcal{T}_h}
\newcommand{\Tht}{\widetilde{\mathcal{T}}_h}
\newcommand{\Thd}{\mathcal{T}_{h,\partial}}
\newcommand{\OT}{\O_{\mathcal{T}}}
\newcommand{\Fh}{\mathcal{F}_h}
\newcommand{\Fgp}{\mathcal{F}_{gp}}

\DeclareMathOperator{\COD}{cod}
\DeclareMathOperator{\TCV}{tcv}
\DeclareMathOperator{\meas}{meas}
\DeclareMathOperator{\dist}{dist}
\DeclareMathOperator{\diver}{div}
\newcommand{\nrm}[1]{\Vert#1\Vert}
\newcommand{\restr}[2]{{\left.\kern-\nulldelimiterspace#1\right|_{#2}}}
\newcommand{\ndef}{n_{\text{def}}}

\begin{document}

\begin{frontmatter}

\title{A coupled high-accuracy phase-field fluid-structure interaction framework for 
Stokes fluid-filled fracture surrounded by an elastic medium}

\author[1,2]{Henry von Wahl}
\ead{henry.von.wahl@uni-jena.de}
\author[3]{Thomas Wick}
\ead{thomas.wick@ifam.uni-hannover.de}

\affiliation[1]{
organization={University of Vienna, Faculty for Mathematics}, 
addressline={Oskar-Morgenstern-Platz 1},
city={Vienna},
postcode={1090},
country={Austria}
}

\affiliation[2]{
organization={Friedrich Schiller University Jena, Institute for Mathematics}, 
addressline={Ernst-Abbe-Platz 2},
city={Jena},
postcode={07743},
country={Germany}
}

\affiliation[3]{
organization={Leibniz University Hannover, Institute for Applied Mathematics}, 
addressline={Welfengarten 1},
city={Hannover},
postcode={30167},
country={Germany}
}

\begin{abstract}
In this work, we couple a high-accuracy phase-field fracture reconstruction
approach iteratively to fluid-structure interaction. The key motivation is
to utilize phase-field modelling to compute the fracture path. A mesh
reconstruction allows a switch from interface-capturing to interface-tracking
in which the coupling conditions can be realized in a highly accurate fashion.
Consequently, inside the fracture, a Stokes flow can be modelled that is
coupled to the surrounding elastic medium. A fully coupled approach is
obtained by iterating between the phase-field and the fluid-structure
interaction model. The resulting algorithm is demonstrated for several
numerical examples of quasi-static brittle fractures. We consider both
stationary and quasi-stationary problems. In the latter, the dynamics arise
through an incrementally increasing given pressure.
\end{abstract}
\end{frontmatter}

\section{Introduction}
\label{sec_intro}
The starting point of the current work is a novel high-accuracy phase-field 
fracture framework with interface reconstructions~\cite{WaWi23_CMAME}. The 
closest related works with similar motivations are~\cite{YNK20, GSF17,
WANG201657}, in which both the phase-field and 
level set/extended finite element (XFEM) approaches are coupled. Next,
\cite{XuYouZhu23} also has the same motivation as we do in which a k-nearest 
neighbor method is used to compute the normal directions of the crack, 
and a subsequent reconstruction. Moreover, 
we mention~\cite{Nguyen2015} and~\cite{LWW17} in which level-set formulations 
from the phase-field function are used to compute normal vectors or construct 
finite element representations of the fracture width. Improvements in terms 
of accuracy of the latter work, namely~\cite{LWW17}, is one specific 
motivation of our current work. Related studies that are also concerned with 
high-accuracy representations of quantities of interest in the fracture 
process zone are~\cite{SARGADO2021113474}, where the accurate resolution 
of stresses in the fracture area was achieved by employing combined 
discretization schemes of finite element and finite volume types,
and~\cite{KUMAR2020104027}, where the phase-field formulation was enhanced
with an additional external driving force.

Various groups have studied fractional phase-field modelling from which
the fracture path is computed~\cite{BourFraMar00, KuMue10, MieWelHof10a, 
MieWelHof10b, BoVeScoHuLa12, AmGeraLoren15, ARRIAGA201833, SARGADO2018458,
WheWiLee20}. Monographs and extended papers of phase-field methods 
for fracture propagation are~\cite{BourFraMar08, AmGeraLoren15, WNN20,
Wic20, HEIDER2021107881, DiLiWiTy22}.
On the other hand, fluid-structure interaction (FSI) in arbitrary 
Lagrangian-Eulerian coordinates (ALE) was proposed in~\cite{HuLiZi81,
DoGiuHa82}; see also the overview~\cite{DoHuePoRo04}. 
Some theoretical work on ALE schemes was done in~\cite{FoNo99,FoNo04}. 
Books and monographs on fluid-structure interaction include~\cite{BuSc06, 
BuSc10, Ri17_fsi, GaRa10, BoGaNe14, BaTaTe13, FoQuaVe09}. The combination 
of phase fields with fluid-structure interaction was subject 
in~\cite{Sun20141, MAO2023111903,Wi16_fsi_pff} and~\cite[Chapter 11.8]{Wic20}.

As previously outlined, we begin with our prior work~\cite{WaWi23_CMAME}, 
in which we used the fracture phase-field to determine the fracture path 
and reconstruct the open fluid-filled crack geometry to have a sharp 
interface between the fracture sub-domain and the surrounding medium. 
This idea is relevant in applications where both a flexible method for 
largely moving or growing interfaces is required, and high accuracy of 
interface conditions is necessary due to different interacting physics 
phenomena. 
One first application can be found in subsurface modelling, for instance,
porous media applications with interfaces~\cite{MiWheWi14, MieheMauthe2015, 
Miehe2015186, LeeWheWi16, YoBou16, HEIDER2018116, CaSaLo18, SaJuaFe17, 
Wilson2016264, AlNoWiWr21_CAMWA, CHUKWUDOZIE2019957, HeiSun20, 
HEIDER2021107881, fei2023phasefield, Heider2016, COSTA2023104040}, 
non-isothermal phase-field fracture~\cite{NoiiWi19, NgSwHeiMa21, NgHeiMa23, 
SUH2021114182},
or computational medicine with fluid-structure interaction~\cite{BADIA20097986,
FoQuaVe09, BoGaNe14, BUKAC2015138, Bukac2015, MR3907987} 
in which interfaces must be resolved with high accuracy in order to keep
discretization errors in the coupling conditions sufficiently small.

The current work's main objective is to extend the one-way coupling proposed
in~\cite{WaWi23_CMAME} to a fully coupled algorithm. However, this includes 
several modifications on both the mathematical and numerical levels. 
Conceptionally, the high-accuracy interface reconstruction framework 
from~\cite{WaWi23_CMAME} is employed in which the fracture path is computed 
with a phase field. Then, the crack opening displacements are computed to
construct the new fluid-filled fracture domain. This allows the application of
an interface-tracking approach with an accurate partition into solid and fluid 
subdomains. In each domain, in principle, any required physics can be posed. 
For simplicity, we focus on elastic solids and linear incompressible (Stokes) 
flow. The overall phase-field fracture fluid-structure interaction is
decoupled and solved iteratively while each sub-problem remains non-linear. 
The phase-field problem consists of displacements and the phase-field 
function that are treated in a fully coupled fashion, where the phase-field 
function is subject to a crack irreversibility condition. These two 
properties result in a non-linear phase-field problem. The fluid-structure 
interaction problem is non-linear due to the ALE approach despite working 
with a linear Stokes flow. In both situations, Newton-type solvers are 
employed. As the reconstruction procedure allows for interface-tracking 
methods, the interface conditions can be prescribed and realized directly on
the interface. This is usually difficult in phase-field methods as it is an 
interface-capturing approach. 

The starting point for our phase-field model is~\cite[Section 2]{MiWheWi13a}, 
in which a porous media mechanics step is proposed using kinematic 
(continuity of principal variables) and dynamic (continuity of fluxes) 
interface conditions between the fracture and the porous medium. In 
\cite[p. 7]{MiWheWi13a}, Gauss' divergence theorem was employed to rewrite the 
sharp interface conditions as domain integrals to account for the property 
in phase-field methods that the interface is given in a smeared fashion rather 
than a precise manifold. This type of approach has been applied in most of 
the references mentioned previously on porous media applications with 
interfaces.

To the best of our knowledge, applying the sharp interface conditions directly
in the phase-field problem, rather than the domain integrals resulting from
applying Gauss' divergence theorem, was never utilized in the published
literature. Our motivation for this approach is as follows. The
original formulation, i.e.,~\cite{FraMar98,BourFraMar00}, starts from a sharp
fracture, which is then approximated for numerical purposes via elliptic
Ambrosio-Tortorelli functionals~\cite{AmTo90,AmTo92}. In this work, go back
to sharp interface approximations. Here, we note that some studies consider
fractures and their flow on dimensionally-reduced objects,
e.g.,~\cite{MaJaRo05,BuYoZu17}, which are, however, given (and fixed) in
advance. In order to develop concepts that can utilize one of the main
intriguing capabilities of phase-field, namely growing, curvilinear
fractures, merging and branching with relative ease, we aim to utilize
phase-field as a state-of-the-art `crack-growth' module. In order to have
all capabilities at hand on the physics in the fracture from first principle
physics laws, e.g., using Stokes or Navier-Stokes, starting with a
dimensionally-reduced fracture model is not an option for us. Our overall aim
is to have the fracture domain and the surrounding medium in the same dimension.
Here, most current works (see the many previous references on phase-field
fracture) employ the physics directly in the phase-field approach, resulting
in multiphysics phase-field fracture. However, the accurate description of
interface conditions is always a point of concern, in addition to relative
coarse sub-meshes in the fracture zone, giving rough finite element
approximations only. Our main motivation is to overcome these last two
limitations by our proposed mesh reconstruction to prescribe interface
conditions on a sharp interface resolved with an interface-tracking method
and, second, to have the flexibility to have fine meshes with sufficiently
small finite element discretization errors. We notice that fluid-structure
interaction with an elastic surrounding medium and Stokes flow in the fracture
is simply a prototype example, and any other physics, e.g., heat transfer,
towards thermo-hydraulic-mechanical fracture, or more complicated fracture
flow models, e.g., Navier-Stokes or non-Newtonian flow, can be utilized
following the same descriptions.

Our first objective is a mathematical formulation of such a pressurized 
phase-field fracture problem statement. This formulation of the phase-field 
problem is then used to couple with the ALE fluid-structure interaction 
problem via our mesh reconstruction approach and the FSI pressure. The 
second objective is the design of a fully coupled algorithm in which we 
iterate between pressurized phase-field fracture and fluid-structure 
interaction that increases the accuracy not only to prescribe the interface 
conditions themselves but also to enhance the accuracy of the overall 
coupled problem to a given tolerance.

As the prescribed pressure drives the phase field, we will couple the 
phase-field problem to the fluid-structure interaction solution through 
the FSI pressure. However, the FSI pressure only exists in the 
reconstructed wide-open fluid domain rather than the entire domain. 
Furthermore, the sharp phase-field crack is a slim domain. Consequently, 
we need to use an averaged fluid-structure interaction pressure in the 
phase-field formulation.

The final algorithm then consist of the computation of the phase-field problem,
the reconstruction of the fluid-filled open crack domain, the computation of
the fluid-structure interaction problem and then the back coupling of an
averaged FSI pressure to the phase-field problem.
This final algorithm is demonstrated for some prototypical scenarios, which
include both stationary and propagating fractures. 
In summary, the novelties are:
\begin{itemize}
  \item Pressurized phase-field fracture problem statement with fracture interface terms;
  \item Design of a fully coupled algorithm for phase-field fluid-structure interaction;
  \item Quasi-static time-dependence with propagating fractures.
\end{itemize}

The outline of this paper is as follows. In \Cref{sec_PFF_FSI}, we present 
our phase-field fluid-structure interaction approach. Next, in 
\Cref{sec.alg}, we discuss the discretization of the phase-field and FSI 
problems, the fluid-filled crack interface reconstruction, the pressure 
coupling and the numerical solution of the discrete problem. This then 
guides us to the design of a final algorithm. In \Cref{sec_tests}, several 
numerical tests substantiate our proposed approach. Finally, our work is 
summarised in \Cref{sec_conclusions}.

\section{A phase-field fluid-structure interaction framework}
\label{sec_PFF_FSI}
In this section, we formulate our high-accuracy phase-field fracture 
fluid-structure interaction framework~\cite{WaWi23_CMAME}. In this framework, 
the fracture shape is determined by a phase-field method. 
In \Cref{sec_PFF}, our phase-field fracture model 
starts from the original work~\cite{BourFraMar00} in which sharp interfaces
are approximated for numerical purposes via elliptic 
Ambrosio-Tortorelli functionals~\cite{AmTo90,AmTo92}. Additional 
physics, e.g., pressure inside the fracture, can be included via two ways, 
namely Gauss' divergence theorem~\cite{MiWheWi13a,MWW19} or keeping 
explicitly the interface integral. The latter is newly investigated in this 
work as it was so far in phase-field fracture due to the smeared zone not being
of interest. We notice that Gauss' divergence theorem is not limited to 
pressure; for instance, temperatures~\cite[Chapter 11]{Wic20}
can be described similarly. The crack irreversibility constraint 
is relaxed using a penalization approach.
We use the phase-field solution to reconstruct the cracked domain, which 
fully resolves the interface between the fluid-filled crack and the intact 
solid domain.
Using this partition, we solve a creeping flow problem inside the crack
coupled to the fluid's elastic medium. In particular, we derive this
fluid-structure interaction problem in \Cref{sec_FSI} by assuming the
surrounding medium as an elastic solid, and the flow in the fracture of
Stokes type.

\subsection{Notation}
In the following, let $\O \subset \RR^d$, $d\in\{2,3\}$, be the total domain
under consideration. Then let $\mathcal{C}\subset\O$ denote the fracture in
our domain and $\SO \subset \O$ denote the intact domain.
In a phase-field approach, the fracture $\mathcal{C}$
is approximated by $\CR\subset\RR^d$ with the help of an elliptic
(Ambrosio-Tortorelli) functional~\cite{AmTo90,AmTo92}.
For fracture formulations posed in a variational setting, this was first
proposed in~\cite{BourFraMar00}. 
The inner fracture boundary is denoted by $\partial\CR$. 
We emphasise that the domains $\SO,\CR$, and the boundary $\partial\CR$ depend
on the choice of the so-called \emph{phase-field regularisation parameter}
$\eps>0$. Details of this parameter are presented below.
Finally, we denote the $L^2(\O)$ scalar product with $(\cdot, \cdot)$.

\subsection{Phase-Field Fracture}
\label{sec_PFF}
We introduce the equations governing the phase-field fracture model used in
this work to model the fluid-filled crack. As we build on our
earlier work~\cite{WaWi23_CMAME}, this is, in most parts, a repetition
of \Cite[Section~2.2]{WaWi23_CMAME}, which we provide for completeness.
However, the phase-field model with sharp interface conditions 
\Cref{form_2_interface_b} is novel in comparison to 
\cite{WaWi23_CMAME} and the literature.

The weak formulations for our phase-field fracture model are given in an 
incremental (i.e., time-discretised) formulation. This is based on a 
pressurised fracture extension of a quasi-static variational fracture model
\cite{FraMar98,BourFraMar00} as presented in~\cite{MiWheWi13a,MWW19}.
We first state a classical formulation and then present the linearised and
regularised formulation we later use in our implementation.

The fracture problem requires two unknowns: the vector-valued displacement
field $\ub$ and the scalar-valued phase-field function $\varphi$. The 
phase-field function indicates the presence of a crack by taking the value
$\varphi=1$ in the intact domain $\SO$, the value $\varphi=0$ inside the
crack $\CR$ and a smooth transition between the two in a region of width
$\eps>0$ around the interface between the open crack and the intact domain, 
denoted now by $\SO$ and $\CR$. In the context of a fluid-filled crack, 
these domains will later be denoted as $\SO$ and $\FL$, respectively.
Since an open crack cannot reseal, the phase field is subject to the
crack irreversibility constraint $\partial_t \varphi \leq 0$.
For the purpose of our phase-field model, the continuous irreversibility 
constraint is approximated through a difference quotient by
\begin{equation*}
  \varphi \leq \varphi^{old}.
\end{equation*}

Later, $\varphi^{old}$ will denote the solution at the previous iteration
step $\varphi^{m-1}$, and the current solution
$\varphi\coloneqq\varphi^m$.

For simplicity, we will assume homogeneous Dirichlet conditions on the
outer boundary $\partial\O$ for the displacement field. Therefore, we
consider the function spaces $W\coloneqq H^1(\O)$,
$\Vb\coloneqq [H^1_0(\O)]^d$ and the convex set
\begin{equation*}
  K\coloneqq \{w\in H^1(\O) |\, w\leq \varphi^{old} \leq 1 \text{ a.e. on }\O\}.
\end{equation*} 
The latter two are the the solution sets, respectively. 
An important object in our work is the interface $\partial\CR$ with the conditions 
(for pressurized fractures):
\begin{align*}
p_{\SO} &= p_{\FL}, \\
\sigma_{\SO}\nb &= \sigma_{\FL}\nb,\\
\end{align*}
where the normal vector $\nb$ points into the crack region. As the interface 
is usually only known as a smeared region in phase-field fracture, 
in the original work~\cite{MiWheWi13a} (see also~\cite[11.4.1.2]{Wic20})
Gauss' divergence theorem was employed yielding for the second condition 
\begin{align}
\int_{\partial\CR} \sigma_{\SO}\nb \, ds
&= \int_{\partial\CR} \sigma_{\FL}\nb \, ds \label{eq_int_1} \\ 
&= \int_{\partial\CR} p\nb\cdot \ub \, ds \label{eq_int_2} \\
&= \int_{\FL} \nabla\cdot (p \ub)\, dx - \int_{\partial\O} p\nb \cdot \ub\, ds \label{eq_int_3}\\
&= \int_{\FL} (\ub\nabla p + p\nabla\cdot \ub)\, dx 
- \int_{\partial\O} p\nb \cdot \ub\, ds\\
&= \int_{\FL} (\ub\nabla p + p\nabla\cdot \ub)\, dx. \label{eq_int_5}
\end{align}
In \eqref{eq_int_2} of this equation chain, we assumed that the leading order of the stress 
$\sigma_{\FL} = -pI + \rho\nu (\nabla v + \nabla v^T)$ is the pressure,
resulting in $\sigma_{\FL} \approx -pI$. In \eqref{eq_int_3}, Gauss' divergence 
theorem is applied. In \eqref{eq_int_5}, homogeneous Dirichlet conditions 
on $\ub = 0$ on $\partial\O$ are assumed.
Despite the fact that the sharp interface representation is the starting point,
in the following we first recapitulate the pressurized phase-field 
fracture model based on \eqref{eq_int_5}, i.e., \Cref{form_1} and \Cref{form_2}, as 
these present our standard model so far. Then, in \Cref{form_2_interface_b} 
the model employing \eqref{eq_int_2} is considered.

To this end, the
phase-field fracture model is given by a coupled variational
inequality system (CVIS)~\cite{Wic20} and reads as follows.
\begin{form}\label{form_1}
Let the pressure $p\in W^{1,\infty}(\O)$, Dirichlet boundary data
$\ub_D$ on $\partial\Omega$, and the initial condition 
$\varphi(0)\coloneqq\varphi_0$ be given.
For the iteration steps $m=1,2,3,\dots,M$, we compute the following problem:
Find $(\ub,\varphi)\coloneqq(\ub^{m},\varphi^{m}) \in \{\ub_D + \Vb\} \times K$,
such that
\begin{subequations}\label{eqn.phase-field.form1}
\begin{align}
 \Bigl(g(\varphi) \;\sigb_s(\ub), \eb( {\wb})\Bigr)
    + ({\varphi}^{2} p, \nabla\cdot  {\wb}) 
+ ({\varphi}^{2} \nabla p, \wb) 
&= 0 \quad \forall \wb\in \Vb,\\
  \begin{multlined}[b]
    (1-\kappa) ({\varphi} \;\sigb_s(\ub):\eb(\ub), \psi {-\varphi}) 
    +  2 ({\varphi}\;  p\; \nabla\cdot  \ub,\psi{-\varphi})
+  2 ({\varphi}\;  \nabla p\; \ub,\psi)\\
    + G_c  \Bigl( -\frac{1}{\eps} (1-\varphi,\psi{-\varphi}) + \eps (\nabla
    \varphi, \nabla (\psi - {\varphi}))   \Bigr)  
  \end{multlined}
    &\geq  0 \quad \forall \psi \in K\cap L^{\infty}(\O).
\end{align}
\end{subequations}
Therein, we have the degradation function
\begin{equation*}
  g(\varphi)\coloneqq (1-\kappa) {\varphi}^2 + \kappa,
\end{equation*}
with the bulk regularisation parameter $\kappa>0$, the phase-field
regularisation parameter $\eps>0$,
the Cauchy stress tensor 
\begin{equation*}
  \sigb_s = 2 \mu \eb(\ub) + \lambda tr(\eb(\ub)) I,
\end{equation*}
with the Lam\'e parameters $\mu,\lambda>0$, the identity matrix 
$I\in\mathbb{R}^{2\times 2}$ and the linearised strain tensor 
\begin{equation*}
  \eb(\ub) = \frac{1}{2} (\nabla \ub + \nabla \ub^T).
\end{equation*}
\end{form}

First, we note that the phase-field regularization parameter $\eps$ is linked 
to the bulk regularization $\kappa$ and the local mesh size parameter $h$; 
see~\cite{Bou99} for a numerical analysis in terms of image segmentation, 
where the main regularization functional (Ambrosio-Tortorelli 
type~\cite{AmTo90,AmTo92}) is similar to phase-field fracture, as well as 
\cite[Section 2.1.2]{BourFraMar00} and~\cite[Section 5.5]{Wic20}.
On the other hand, $\eps$ is known as a length scale in engineering and 
can be interpreted as a fixed number representing micro-cracks around 
the fracture zone; see, e.g.,~\cite{MieWelHof10a,Nguyen16}. 
For typical numerical-computational convergence studies, a decision 
must be made about which type of convergence is to be studied. For 
discretization errors only, $\eps$ must be fixed. If the interaction of 
$\eps$ and $h$ towards sharp interface limits is of interest,
$\eps$ should also vary. However, significant computational resources, mesh 
adaptivity and parallel computing are needed to study asymptotic 
limits~\cite{HeiWi18_pamm}. Furthermore, $\eps$ impacts crack nucleation, 
e.g.,~\cite{Bour07,TANNE201880}, which limits its choice. Finally, we note 
that as an extension, $\eps$ can be made arbitrary at the cost of introducing 
additional parameters in the total energy-functional~\cite{SARGADO2018458}.

Next, we note that the system \eqref{eqn.phase-field.form1} does not explicitly
contain time-derivatives. In our setting, the time $t$ enters through 
time-dependent, i.e., incrementally dependent since we work in a quasi-static 
regime, right-hand side forces, e.g., a 
time-dependent pressure force $p\coloneqq p(t)$.
Due to the quasi-static nature of this problem formulation, we derive 
the time-discretised formulation with some further approximations. 
Our first approximation relaxes the non-linear behaviour in the
first term $g(\varphi) \;\sigb_s(\ub)$ of the displacement equation by using 
\begin{equation*}
\varphi \approx \varphi^{m-1}, 
\end{equation*}
yielding $g(\varphi^{m-1})$.
This idea is based on the 
extrapolation introduced in~\cite{HWW15} and is numerically justified 
specifically for slowly growing fractures~\cite[Chapter 6]{Wic20}, 
while counter-examples for fast-growing fractures were found in~\cite{Wic17}.

As $\varphi^{m-1}$ introduces an approximation error of order 
$O(\tau_m - \tau_{m-1})$ for smaller steps, the error vanishes, or 
fully monolithic schemes (again~\cite{Wic17}) or an additional 
iteration~\cite{KMW23} must be introduced.

The second approximation is related to the inequality constraint. In this 
work, we relax the constraint by simple penalisation~\cite{MWT15} (see also 
\cite[Chapter 5]{Wic20}), i.e., 
\begin{equation*}
\varphi \leq \varphi^{m-1} \quad\rightarrow\quad \gamma(\varphi - \varphi^{m-1})^+.
\end{equation*}
Here, $(x)^+ = x$ for $x>0$ and $(x)^+ = 0$ for $x\leq 0$, 
and where $\gamma>0$ is a penalty parameter.
We then arrive at the regularised scheme
\begin{form}\label{form_2}
Let $p\in W^{1,\infty}(\O)$ and the initial condition $\varphi(0)\coloneqq\varphi_0$ be given.
For the iteration steps $m=1,2,3,\dots,M$, we compute:
Find $(\ub,\varphi)\coloneqq(\ub^m,\varphi^m) \in \Vb \times W$ such that
\begin{subequations}
\begin{align}
  \Bigl(g(\varphi^{m-1})\;\sigma(\ub), \eb({\wb})\Bigr)
    +({\varphi^{m-1}}^{2} p, \nabla\cdot{\wb}) 
    + ({\varphi^{m-1}}^{2} \nabla p, \wb) 
    &= 0
  \quad \forall \wb\in \Vb,\\
  \begin{multlined}[b]
    (1-\kappa) ({\varphi} \;\sigma(\ub):\eb(\ub), \psi) 
      +  2 ({\varphi}\;  p\; \nabla\cdot  \ub,\psi)
      +  2 ({\varphi}\;  \nabla p\; \ub,\psi)\\
    +  G_c  \Bigl( -\frac{1}{\eps} (1-\varphi,\psi) 
      + \eps (\nabla\varphi, \nabla\psi) \Bigr) 
      + (\gamma(\varphi - \varphi^{m-1})^+,\psi)
  \end{multlined}
   &=0\quad\forall \psi\in W.
\end{align}
\end{subequations}
Note that we assume a varying pressure $p$, which requires 
to use the pressure gradient $\nabla p$~\cite[(2.30) and (2.31)]{MiWheWi14}
in extension to our previous work~\cite{WaWi23_CMAME}.
\end{form}

\begin{remark}
We call the steps $\tau_m$ with index $m$ in this work `iteration' steps, 
which are used to satisfy the irreversibility constraint. In previous works, 
they were introduced as incremental steps, pseudo-time steps,  or time steps. 
We choose iteration steps in this work to better distinguish from (time) steps 
$t_n$ with index $n$ used to advance the overall coupled system of phase-field 
mesh reconstruction and fluid-structure interaction.
\end{remark}

The previous formulation is based on an interface law is formulated 
as domain integral using Gauss' divergence theorem
\cite[Section 2]{MiWheWi13a}. An equivalent formulation can be given by 
directly using the interface term. This then is given as follows.
\begin{form}\label{form_2_interface_b}
Let $p\in L^{\infty}(\O)$ and the initial condition
$\varphi(0)\coloneqq\varphi_0$ be given. For the iteration steps 
$m=1,2,3,\dots,M$, we compute:
Find $(\ub,\varphi)\coloneqq(\ub^m,\varphi^m) \in \Vb \times W$ such that
\begin{subequations}\label{eqn.form_2_interface_b}
\begin{align}
  \Bigl(g(\varphi^{m-1})\;\sigma(\ub), \eb({\wb})\Bigr)
    + \int_{\partial\CR} {\varphi^{m-1}}^2 p \nb \cdot \wb \, ds &= 0
  \quad \forall \wb\in \Vb,\label{eqn.form_2_interface_b.a}\\
  \begin{multlined}[b]
    (1-\kappa) ({\varphi} \;\sigma(\ub):\eb(\ub), \psi)
      +  2 \int_{\partial\CR} \varphi p \nb \cdot \ub \psi \, ds\\
    +  G_c  \Bigl( -\frac{1}{\eps} (1-\varphi,\psi) 
      + \eps (\nabla\varphi, \nabla\psi) \Bigr) 
      + (\gamma(\varphi - \varphi^{m-1})^+,\psi)
  \end{multlined}
   &=0\quad\forall \psi\in W.\label{eqn.form_2_interface_b.b}
\end{align}
\end{subequations}
Here, the normal vector $\nb$ points into the crack region. 
\end{form}
To the best of our knowledge, \Cref{form_2_interface_b} has not been used 
in the literature thus far. This can be attributed to the fact that it 
contradicts the concept of a phase field in the sense that a phase-field 
model should not require explicit knowledge of the interface position. 
However, in our proposed coupling with fluid-structure interaction, and 
with the aim to provide a high-accuracy framework, \Cref{form_2_interface_b} 
is of interest. Consequently, we will use \eqref{eqn.form_2_interface_b}  
as the system of equations to compute the phase-field fracture in each 
iteration. As we are in a quasi-stationary regime here, we compute
each phase-field with a total of $M=5$ iteration steps. In each of these
iteration steps, we then use Newton's method to solve the resulting non-linear
problem up to a tolerance of $10^{-8}$.
Time dependence will only enter our system 
through time-dependent pressure and consequently changing crack domains.

In~\cite{WaWi23_CMAME}, we presented a number of methods to reconstruct the 
geometry of an open, fluid-filled crack resulting from a phase-field model. 
Here, we will focus on the approach where the resulting computational mesh 
resolves this geometry. With the mesh containing the resolved interface 
between the open crack and the surrounding solid at hand, we can now describe 
the fluid-structure interaction problem between the fluid-filled crack and 
the solid domain surrounding it. With our explicit interface reconstructions 
of the fracture surface, the flow problem is coupled via an interface-tracking 
approach to the surrounding solid, and we arrive at a classical 
fluid-structure interaction model. In order to couple flow and solids, we 
discuss the arbitrary Lagrangian-Eulerian approach below.

\subsection{Stationary Fluid-Structure Interaction}
\label{sec_FSI}
In this sub-section, we model the previously announced fluid-structure
interaction (FSI) problem in arbitrary Lagrangian-Eulerian (ALE) coordinates 
using variational monolithic coupling in a reference configuration
\cite{HrTu06a,Du07,Wi11_phd,Ri17_fsi}.
We assume a creeping flow inside the fluid-filled crack. Consequently,
we consider the Stokes equations for the fluid model.
Nevertheless, our method is not limited to Stokes and can deal with 
Navier-Stokes in the same fashion. However, from a physics point of view,
having quasi-stationary settings in mind in this paper, we assume 
a creeping flow inside the crack and thus the Stokes model is reasonable
to be employed.

For the FSI problem, we consider a domain $\O\subset\RR^d$ divided into
a $d$-dimensional fluid domain $\FL$, a $d$-dimensional solid domain
$\SO$ and a $d-1$-dimensional interface $\IN$ between the two, such that
$\O=\FL\dot\cup\IN\dot\cup\SO$. Furthermore, let $\hat{\O}, \FLref,
\SOref$ and $\INref$ be the corresponding domains in a reference
configuration. Similarly, $\vbref, \pref, \ubref$ and $\xbref$ denote the
velocity, pressure, deformation and coordinates in the reference configuration.
In our setting of a fluid-filled crack, the fluid domain
is the interior of the crack $\FL=\CR$, the solid is the intact medium
$\SO$ and the interface is the crack boundary $\IN=\partial\CR$.
Furthermore, for a domain $X$, we denote the $L^2(X)$ scalar product
with $(\cdot,\cdot)_X$.

Using the reference domains $\SOref$ and $\FLref$ leads to
the well-established formulation in ALE coordinates~\cite{HuLiZi81,DoGiuHa82}. 
To obtain a monolithic formulation we need to specify the transformation
$\hcalA_f$ from the reference configuration to the physical domain in
the fluid-domain. On the interface $\INref$ this transformation
is given by the structure displacement:
\begin{equation*}
  \hcalA_f(\xbref,t)\big|_{\INref} = \xbref+\ubref_s(\xbref,t)\big|_{\INref}. 
\end{equation*}
On the outer boundary of the fluid domain
$\partial\FLref\setminus\INref$ it holds $\hcalA_f=\text{id}$.
Inside $\FLref$, the transformation should be as smooth and regular as
possible, but apart from that it is arbitrary. Thus we harmonically
extend $\ubref_s|_{\SOref}$ to the fluid domain $\FLref$ and
define $\hcalA_f\coloneqq \text{id}+\ubref$ on $\FLref$,
where $id(\xbref) = \xbref$ in $\hcalA_f\coloneqq \text{id}+\ubref
\coloneqq \hcalA_f(\xbref,t)\coloneqq \text{id}(\xbref)+\ubref(\xbref,t)$
such that
\begin{equation*}
  (\hat \nabla \ubref_f,\hat \nabla\psiref)_{\FLref} = 0,\quad
  \ubref_f=\ubref_s\text{ on }\INref,\quad
  \ubref_f=0\text{ on }\partial\FLref\setminus\INref.
\end{equation*}
Consequently, we define a continuous variable $\ubref$ on  all $\Omega$
defining the deformation in $\SOref$ and supporting the
transformation in $\FLref$. By skipping the subscripts and
since the definition of $\hcalA_f$ coincides with the definition of
the solid transformation $\hcalA_s$, we define on
all $\Oref$:
\begin{equation*}
  \hcalA(\xbref, t)\coloneqq \xbref + \ubref(\xbref, t),\quad
  \hat F(\xbref, t) \coloneqq \hat\nabla \hcalA= I+\hat\nabla \ubref(\xbref, t),\quad
  \hat J\coloneqq \text{det}(\hat F).
\end{equation*}

With this at hand, the weak formulation of the stationary fluid-structure
interaction problem~\cite{RiWi10} is given as follows.
\begin{form}[Stationary fluid-structure interaction]
\label{fsi:ale:stationary}
Let $\Vbref$ be the subspace of $\bm{H}^1(\Oref)$ with trace zero on
$\hat\Gamma^D\coloneqq \hat\Gamma_f^D\cup\hat\Gamma_s^D$ and 
$\hat L\coloneqq L^2(\Oref)/\mathbb{R}$. Find $\vbref\in \Vbref$, $\ubref\in
\Vbref$ and $\hat p\in \hat L$, such
\begin{subequations}\label{eqn.ale-fsi}
  \begin{align}
    (\hat J\sigbref_f\hat F^{-T},\hat \nabla\phiref)_{\FLref}
      + (\hat J\sigbref_s\hat F^{-T},\hat \nabla\phiref)_{\SOref}
      &= (\rho_f \hat J\fbref,\phiref)_{\FLref} &&\forall\phiref\in\Vbref,\\
    - (\vbref,\psiref)_{\SOref} + 
      (\alpha_u \hat \nabla \ubref,\hat \nabla\psiref)_{\FLref}
      &=0 &&\forall\psiref\in\Vbref,\\
    (\widehat{\diver}\,(\hat J\hat F^{-1}\vbref_f),\xiref)_{\FLref} 
      &=0&&\forall\xiref\in \hat L,
  \end{align}
\end{subequations}  
with a right-hand side fluid force $\fbref\in L^2(\FLref)$ and the 
harmonic mesh extension parameter $\alpha_u>0$. 
Finally, the Cauchy stress tensor in the solid is
defined in \Cref{form_1} and we use
$\hat J\sigbref_s\hat F^{-T} \coloneqq  \sigb_s$.
The ALE fluid Cauchy stress tensor $\sigbref_f$ is given by 
\begin{equation*}
  \sigbref_f \coloneqq  -\hat p_fI +\rho_f\nu_f(\hat\nabla \vbref_f \hat F^{-1}
  + \hat F^{-T}\hat\nabla \vbref_f^T),
\end{equation*}
with the kinematic viscosity $\nu_f>0$ and the fluid's density $\rho_f>0$.
\end{form}

\section{Interface reconstruction, and coupled algorithm, and discretisation}
\label{sec.alg}

The coupled phase-field fluid-structure interaction problem is discretised with 
Galerkin finite elements and is non-linear. These non-linearities arise 
through the overall coupling, the crack irreversibility condition, and the 
ALE transformation. The two governing problems, namely 
phase-field and fluid-structure interaction, are solved in an iterative
fashion, in which both remain non-linear in their nature. Here, Newton-type
solvers are employed.

\subsection{Interface reconstruction}
\label{sec.alg:subsec.interface}
Following our work in \Cite{WaWi23_CMAME}, the open crack domain is 
reconstructed based on the crack opening displacements or aperture of the 
crack. This can be computed by
\begin{equation*}
 \COD(\xb) = \llbracket \ub\cdot\nb \rrbracket
 \simeq \int_{\ell^{\xb,\bm{v}}} \ub(\xb) \cdot \nabla \varphi(\xb) \dif s,
\end{equation*}
where $\ell^{\xb,\bm{v}}$ is a line through $\xb$ along the vector $\bm{v}$ 
\cite{CHUKWUDOZIE2019957}. To avoid unphysical behaviour, due to large 
relative errors near the crack tips, where the exact COD is zero, we disregard
COD values below $h/10$. Having computed the COD at a number of points, and 
with knowledge of the crack centreline, we have a given number of points on 
the crack interface. These points are then connected into a curve to describe 
the geometry of the open crack geometry, which we refer to below as 
$\Omega_\text{cod}$. This is then meshed and taken as the reference domain 
for the fluid-structure interaction solver.

\subsection{Pressure coupling}
The phase-field dynamics are primarily driven by the pressure $p$, as can be 
seen in \Cref{form_2} and \Cref{form_2_interface_b}. Therefore, to couple the 
phase-field and fluid-structure interaction problems, we need to be able to 
use the FSI pressure ($\hat p_f$ in \Cref{fsi:ale:stationary}) in the 
phase-field problem \eqref{eqn.form_2_interface_b}. However, the fluid 
pressure lives in the open crack fluid domain $\FL$, and the crack pressure 
needs to be in the thin crack $\CR$. Furthermore, while the fluid domain has 
a width $\pazocal{O}(1)$, the crack width is $\pazocal{O}(h)$. The latter is 
due to the construction of the basic phase-field model, see e.g., 
\cite{BourFraMar00}, by using Ambrosio-Tortorelli functionals 
\cite{AmTo90,AmTo92} to approximate the fracture area. On the other hand, 
working explicitly with the crack fluid domain $\FL$ with a sharp interface, 
the main interest is to get away from a $h$ dependence of the domain 
approximation.

We, therefore, construct a vertically averaged pressure by computing the mean 
of the pressure on $\pazocal{O}(1/h)$ equally distant lines normal to the
crack centreline and interpolating this into a finite element 
function on the phase-field mesh. Using unfitted finite element
quadrature~\cite{BCH14}, we can integrate on arbitrary lines not resolved 
by the mesh. The resulting pressure then models the effects of 
the fluid-structure interaction pressure inside the thin crack domain $\CR$.

\subsection{Final algorithm}
\label{sec.alg:subsec:final-alg}
In this key sub-section, we propose our fully-coupled scheme for 
combining and coupling the phase-field problem, domain reconstruction and
fluid-structure interaction problems. This leads us to the following algorithm:

\begin{algorithm}[H]
  \KwData{Initial domain $\Omega$ with crack configuration, and a time-dependent background pressure $\overline{p}(t)$.}
  Compute initial phase-field \Cref{form_2_interface_b} in $\Omega$\;
  Initialize $p^\ast=0$\;
  \For{Time point indices $n = 0,\dots, N$}{
    $t^n = n \Delta t$\;
    \For{$k = 1,\dots, K$}{
      Compute phase-field \Cref{form_2_interface_b} using the modified pressure $p = \overline{p}(t^n) + p^*$\;
      Reconstruct the fluid filled crack geometry $\Omega_\text{cod}$ from
      $\Omega$ using the COD, resulting in the current physical domain with a 
      resolved fluid and solid partition\;
      Solve the Stokes FSI \Cref{fsi:ale:stationary} in $\widehat{\Omega}\coloneqq
      \Omega_\text{cod}$\;
      Compute averaged pressure $p^*$ along lines perpendicular to crack\;
      }}
  \caption{Full time stepping loop}
  \label{alg.full-loop}
\end{algorithm}

In line 8 of \Cref{alg.full-loop}, we observe that we use a combined pressure 
consisting of an external, hydrostatic background pressure and the pressure 
resulting from the fluid-structure interaction problem in Line 6. This is 
because the pressure in the Stokes in the Stokes FSI problem is only defined 
uniquely up to a constant (as it is well known, see, e.g.,~\cite{Temam2001}), 
and this constant is fixed by prescribing a mean value of zero. As this mean 
value is essentially arbitrary, it is valid to change this constant in this 
setting. Furthermore, since \Cref{form_2_interface_b} requires an 
explicit description of $\CR$, we need to construct a new mesh for the 
phase field in line 8 to account for the moving crack tip.

Finally, we comment more explicitly on all iteration loops 
in \Cref{alg.full-loop}. In line 3, with index $n$ and final index $N$, 
we denote the incremental (time-like) procedure. In lines 5-10, the overall
phase-field mesh reconstruction fluid-structure interaction system is solved 
$K$ times indexed by $k$. We carry out $K=5$ iterations.
In line 6, the phase-field system is solved at each $t_n, n=1,\ldots,N$
for iteration steps $m=1,\ldots,M$ (at each $t_n$) in order to satisfy the 
crack irreversibility constraint. At each Index $m$, \Cref{form_2_interface_b} 
is non-linear and solved by Newton's method (indexed by $j$, irrelevant here)
up to a tolerance $10^{-8}$. Similarly, in line 8, the non-linear 
fluid-structure interaction system \Cref{fsi:ale:stationary} is solved 
by Newton's method (iteration index $j$, irrelevant here) up to a tolerance 
$10^{-8}$. In both cases, direct sparse linear solvers are employed;
see \Cref{sec_num_sol}.

\subsection{Spatial discretization}
We consider shape-regular simplicial meshes of the phase-field and 
fluid-structure interaction domains with characteristic mesh size $h$. In the 
case of the phase-field problem, the crack then consists of a slit of width 
$2h$ and specified length, and the mesh is refined towards the crack, such 
that $h_\text{max} = 100 h_\text{crack}$. On this mesh, we consider continuous, 
piecewise-linear finite elements $\PP^1$ for both the phase-field and 
deformation variables.

Similarly, we mesh the fluid-structure interaction geometry with a local mesh 
parameter for the fluid region (the open crack) of $h=h_\text{crack}$ and the 
surrounding intact solid with $h=h_\text{max}$. On this mesh, we consider the 
lowest order Taylor-Hood elements: vector-valued $\PP^2$ elements for the 
velocity and deformation and $\PP^1$ elements for the pressure.

\subsection{Numerical solution}
\label{sec_num_sol}
The non-linear systems resulting from the discretisation of the phase-field 
problem in \eqref{eqn.form_2_interface_b} and the fluid-structure interaction 
problem \eqref{eqn.ale-fsi} are solved using Newton's scheme. The Jacobian is 
computed through automatic differentiation by \texttt{NGSolve}\footnote{For 
further details on \texttt{NGSolve}'s implementation of Newton scheme we refer
to the documentation \url{https://docu.ngsolve.org/latest/index.html}.}. The 
resulting linear systems are solved using the direct solver \texttt{pardiso} 
through the Intel MKL library. We refer to the available 
code~\cite{vWW23_zenodo} used to compute the numerical results below for 
further details.

\section{Numerical tests}
\label{sec_tests}
We present a number of numerical examples of our final algorithm presented 
in \Cref{sec.alg}. Our examples are implemented using 
\texttt{Netgen/NGSolve}~\cite{Sch97,Sch14} with the addition add-on 
\texttt{ngsxfem}~\Cite{LHPvW21} for unfitted finite elements. We use the 
unfitted finite element technology to compute the COD and vertically averaged
pressure  along arbitrary lines, which do not need to be resolved by the mesh.
The code  implementing these examples is freely available; see also the Data
Availability Statement below.

\subsection{Example 1: Sneddon test with Formulation~\ref*{form_2_interface_b}}
As discussed above, \Cref{form_2_interface_b} of the pressurised crack 
propagation problem is necessary in order for us to couple the FSI pressure, 
which is only available inside the crack, to the phase-field model. Since 
this formulation of the problem has not been previously used in the literature
in practice (c.f. our discussion in \Cref{sec_intro}), we start with 
Sneddon's test~\cite{Sne46, SneddLow69} to convince ourselves of the
validity of this formulation.

\subsubsection{Configuration}
The domain of interest is $\O=(-2, 2)^2$, and because the data driving the 
crack is constant, the problem is stationary. The mechanical parameters are 
Young's modulus and Poisson's ratio, which we take as $E_s = 10^5$ and 
$\nu_s = 0.35$. The relation to the Lamé parameters is 
$\mu_s = E_s / (2(1 + \nu_s)$ and 
$\lambda_s = \nu_s E / ((1 + \nu_s)( 1 - 2 \nu_s))$. The applied pressure is 
$p = 4.5 \times 10^3$, and we choose the critical energy release rate as 
$G_c = 500$. The boundary conditions applied are
\begin{equation*}
  \ub = 0\quad\text{and}\quad \eps\dn \varphi = 0 \quad\text{on }\partial\O.
\end{equation*}
The initial condition for the phase-field is computed by setting $\phi(\xb)=1$ in $\CR_0$ and zero everywhere else. In order for the resulting function to satisfy the governing equation, we solve 
\eqref{eqn.form_2_interface_b.b} for $M$ iterations with $\ub=0$. This is to avoid an artificial increase of toughness of the crack at it's tips, as seen in~\cite{SARGADO2018458}.

The initial crack is given by $\CR_0 = (-l_0, l_0)\times(-h, h)$ with the 
crack length $2l_0=0.4$ and the local mesh size $h = h_\text{crack}$.
The mesh of the domain is non-uniform and constructed such that elements
in $\O\setminus\CR$ have a diameter of up to $h_\text{max} = 100h_\text{crack}$.

The phase-field penalisation parameter is chosen as 
$\gamma = 100 \cdot h^{-2}$, and the phase-field regularisations parameters 
are $\kappa = 10^{-10}$ and $\eps = 0.5\cdot h^{1/2}$. We iterate the 
phase-field problem for five iteration steps, i.e., $M=5$, to arrive at the
stationary solution but keep $\partial\CR$ fixed. Furthermore, we fix 
$\varphi, \varphi^{n-1}\equiv 1$ in the boundary integrals in 
\eqref{eqn.form_2_interface_b}, such that these terms can be ignored. 
The finite element spaces consist of $H^1$-conforming, 
piecewise linear functions for both the displacement and the phase field.

\subsubsection{Results}
\begin{figure}
  \centering
  \includegraphics{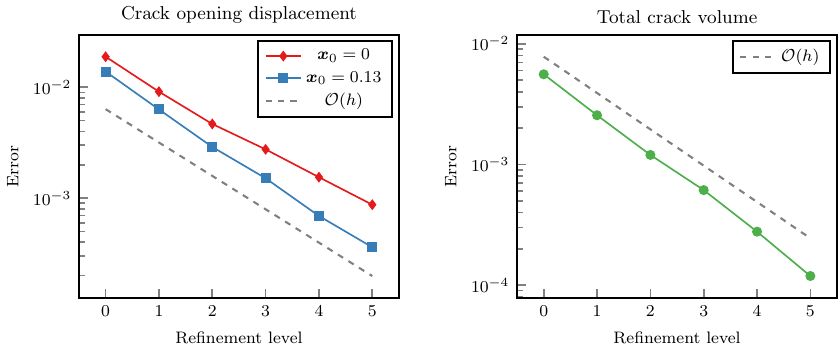}
  \caption{Crack opening displacement convergence for Sneddon's test computed
    computed using \Cref{form_2_interface_b}.}
  \label{fig.results.sneddon}
\end{figure}

We compute the problem over a series of meshes constructed with 
$h_\text{crack}=h_0 \cdot 2^{-\ell}$, $\ell=0,\dots, 5$ and $h_0 = 0.02$. 
As quantities of interest, we consider the crack opening displacements at 
$\xb_0 = 0, 0.13$ and the total crack volume, for which we have the
asymptotically valid analytical expression
\begin{equation*}\COD(x)=4\frac{(1-\nu_s^2)l_0 p}{E_s}\biggl(1-\frac{x^2}{l_0^2}\biggr)^{1/2}
  \quad\text{and}\quad
  \TCV = 2\pi \frac{(1 - \nu_s^2)l_0^2p}{E_s},
\end{equation*}
respectively. The resulting errors can be seen in \Cref{fig.results.sneddon}.

\subsection{Example 2: Modified Sneddon test with Algorithm~\ref*{alg.full-loop}}
We consider a set-up inspired by the previous example, but including 
\Cref{alg.full-loop} in a stationary setting. That is, we iterate between 
the phase-field and FSI problems, thereby coupling the stationary phase-field 
problem to the pressure resulting from the FSI problem.

\subsubsection{Configuration}
The spatial configuration is as in the previous example. The stationary Stokes 
FSI \Cref{fsi:ale:stationary} is driven by the right-hand side forcing term
\begin{equation*}
  \fbref = c_1 \exp(-c_2 ((x - x_0)^2 + (y - y_0)^2))\bm{e}_x,
\end{equation*}
with $c_1=0.02, c_2=1000$ and $(x_0, y_0)=(l_0 / 4, 2 a / 3)$ where $l_0=0.2$ 
is the half length of the crack and $a>0$ is the half-height of the resulting 
ellipse shaped crack, i.e., $\COD(0) / 2$. This results in a negative pressure 
at the left tip and a positive pressure at the right crack tip. Consequently, 
we expect the crack resulting from the coupled computation to be smaller on 
the left end of the crack and wider towards the right of the crack. The 
material parameters for the fluid are $\nu_f=0.1$ and $\rho_f=10^3$. Finally, 
we choose the harmonic mesh extension parameter as $\alpha_u=10^{-14}$.

In order to see the effects of the averaged FSI pressure $p^*$ on the 
phase-field problem, we chose the data as $E=10^2$ and $p=4.5$. Consequently, 
the resulting crack shape is identical to the previous example if the FSI 
pressure is ignored.

\subsubsection{Results}

\begin{figure}
  \centering
  \begin{minipage}[b]{0.49\textwidth}
    \centering
    \includegraphics{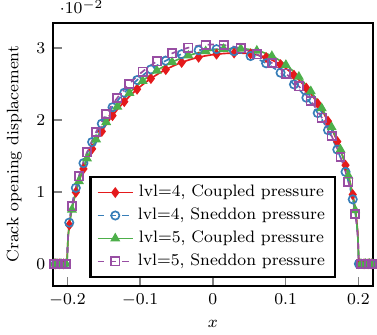}
  \end{minipage}
  \begin{minipage}[b]{0.49\textwidth}
    \centering
    \includegraphics[height=4.5cm, trim=150 0 150 0, clip]{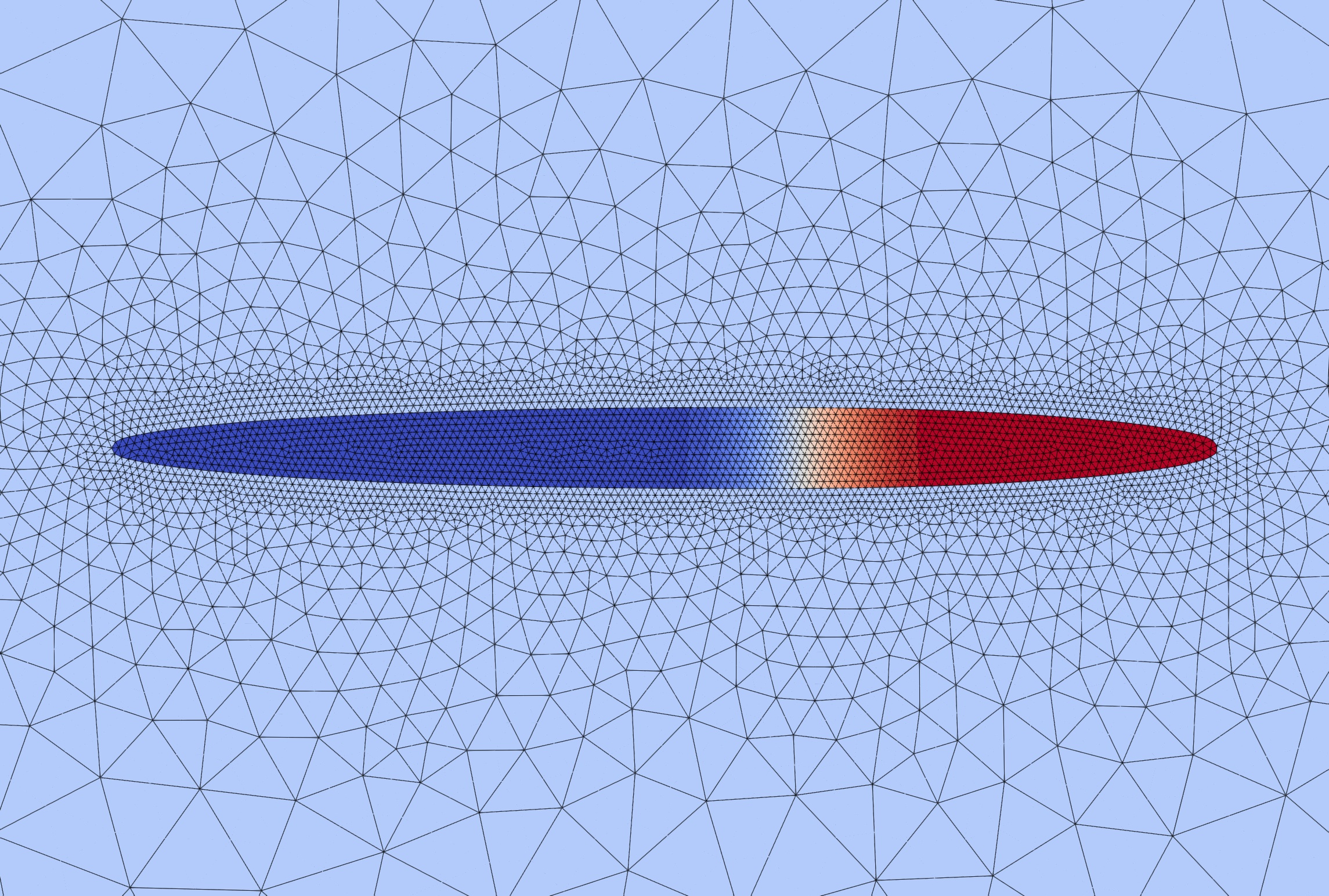}
    \includegraphics[height=4.5cm]{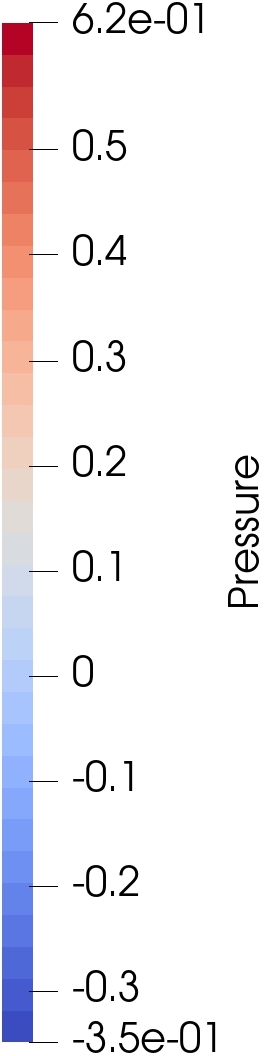}
    \vspace*{23pt}
  \end{minipage}
  \caption{Example 2: Fluid-filled phase-field crack coupled with the 
    fluid-structure-interaction pressure. Left: Resulting crack shape based 
    on the crack opening displacements after four sub-iterations between the 
    phase-field and the fluid-structure-interaction problem with 
    $lvl\in\{4, 5\}$ levels of mesh refinement. Right: Fluid-Structure 
    interaction pressure on the reconstructed domain used to drive the 
    divergence from Sneddon's test on mesh level three.}
  \label{fig.results.ex2:shape-pressure}
\end{figure}

The resulting crack shape on the two finest meshes (four and five 
levels of mesh refinement, respectively) from the coupled FSI 
pressure and Sneddon background pressure can be seen on the left of 
\Cref{fig.results.ex2:shape-pressure} (solid lines) together with the 
resulting Sneddon crack, i.e., no FSI pressure (dashed lines). On the right 
of \Cref{fig.results.ex2:shape-pressure}, we show the FSI pressure inside the 
reconstructed crack. We observe that where the FSI pressure is positive, the 
resulting crack is larger than the Sneddon crack, and where it is negative, 
the COD is smaller, as is to be expected.

\begin{figure}
  \centering
  \includegraphics{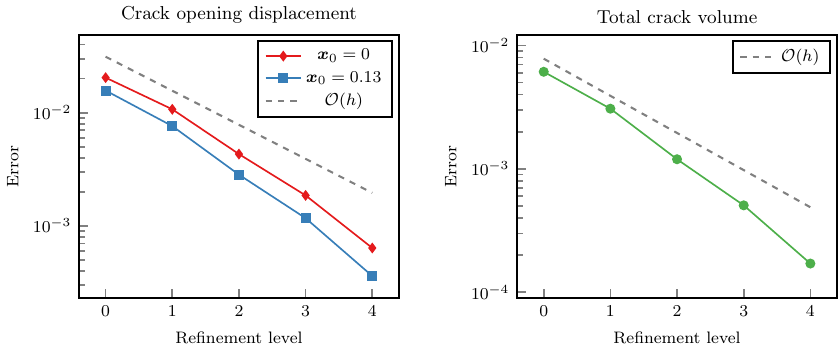}
  \caption{Example 2: Fluid-filled phase-field crack coupled with the 
    fluid-structure-interaction pressure. Left: Numerical convergence for 
    the crack opening displacements in two points. Right: Numerical 
    convergence of the total crack volume.}
  \label{fig.results.ex2:convergence}
\end{figure}

In \Cref{fig.results.ex2:convergence}, we see the convergence of the crack 
opening displacements in the points $x=0$ and $x=0.13$, as well as the 
convergence of the total crack volume towards the reference value. The 
reference value is the result computed on the fines mesh (five levels of 
mesh refinement). We again observe approximately linear mesh convergence 
for all three quantities.

\subsection{Example 3: Crack starting from a boundary with fluid flow boundary conditions}
In this example, we consider a crack from the boundary of the domain, rather than a crack fully contained in the domain. This contact of the crack domain with the boundary adds technical challenges to the reconstruction procedure. 

\subsubsection{Configuration}
The domain of interest is $\Omega=(-1, 1)^2$. The initial crack is given by $\CR_0 = (-1, 0)\times(-h, h)$.
The mechanical parameters are $E_s = 5\cdot 10^4$ and $\nu_s = 0.35$. We choose the model parameters in the phase-field problem as $G_c=500$, $\kappa=10^{-10}$, $\gamma = 100 h^{-2}$, and $\varepsilon= 2 h$. The phase-field background pressure is $10^3$.

In contrast to the previous examples, the crack touches the outer boundary. Therefore, we have to proceed non-homogeneous Dirichlet boundary conditions for the deformation field. Following~\cite[Section~6.3]{Wic20} these are given as $\ub_x = 0$, and
\begin{equation*}
  \ub_y = \begin{cases}
    c \sqrt{\sqrt{x^2 + y^2} - x} &\text{if }y>0\\
    -c \sqrt{\sqrt{x^2 + y^2} - x} &\text{if }y\leq0,
  \end{cases}
\end{equation*}\
with $c=0.08/\sqrt{2}$. For the phase-field, we again have homogenous Neumann boundary conditions.

As the crack touches the outer boundary, we can also consider a different set-up for the fluid-structure interaction problem. Rather than specifying a volume force, we specify an inflow boundary condition on the outer boundary of the open crack. Specifically, we set
\begin{equation*}
  \vbref_{f} = \begin{pmatrix}
    v_\text{max} (d / 2 - y)(d / 2 + y) / d^2\\
    0
  \end{pmatrix},
\end{equation*}
where $v_\text{max}=0.001$ is the maximal inflow speed and $d=\COD(x=-1)$ is the size of the crack opening at the boundary. For the FSI deformation, we set homogeneous Dirichlet Conditions on the top, bottom and inflow boundaries, and homogeneous Neumann boundary conditions on the left and right boundaries. The body forcing term is set as $\fbref=0$.

\subsubsection{Results}
\begin{figure}
  \centering
  \begin{minipage}[t]{.49\textwidth}
    \centering
    \includegraphics[height=5cm]{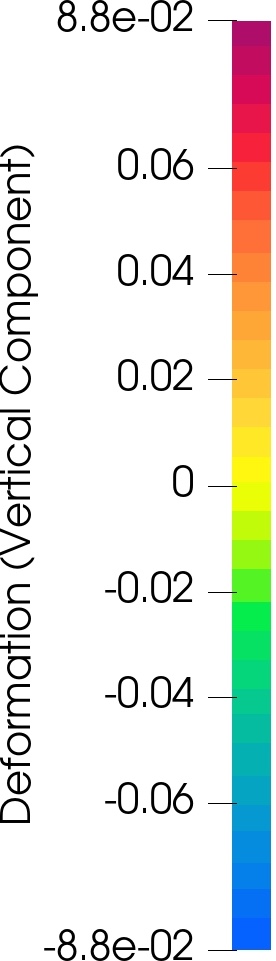}
    \hspace*{5pt}
    \includegraphics[height=5cm]{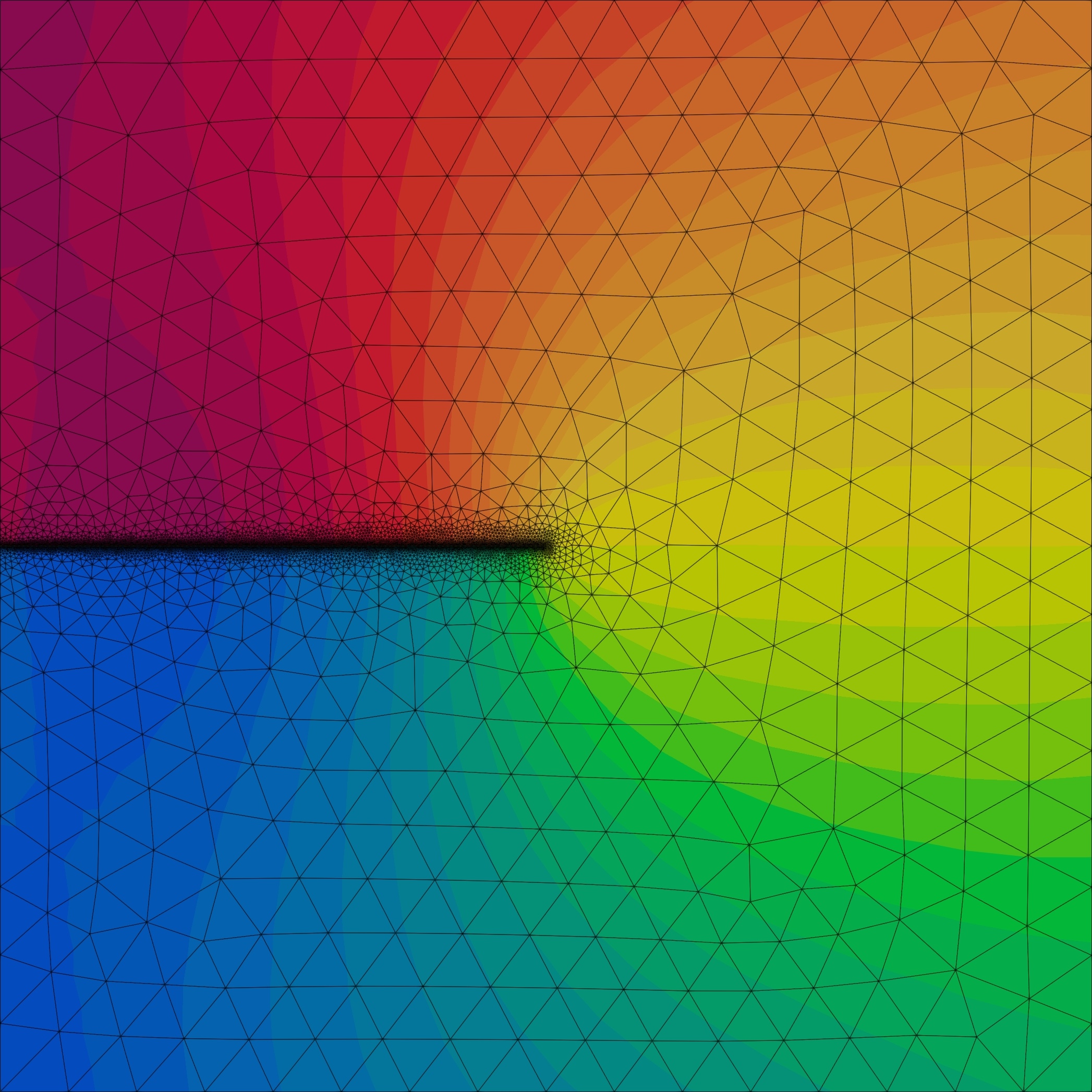}
  \end{minipage}
  \begin{minipage}[t]{.49\textwidth}
    \centering
    \includegraphics[height=5cm]{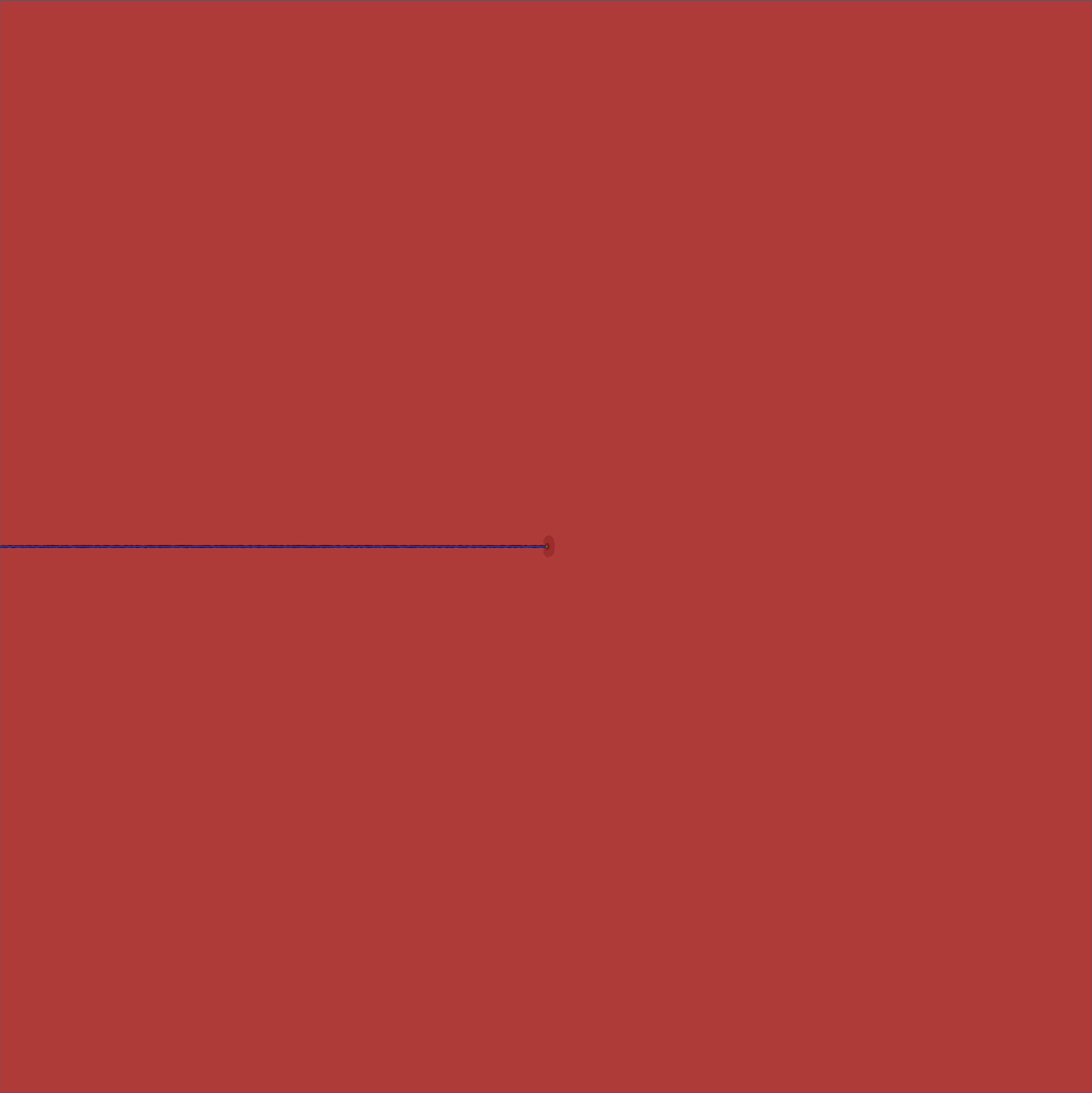}
\end{minipage}
  \caption{Example 3: Phase-field-deformation (vertical component) and phase-field after three levels of mesh refinement.}
  \label{fig.results.ex3:phase-field}
\end{figure}

\begin{figure}
  \centering
  \begin{minipage}[t]{.32\textwidth}
    \centering
    \includegraphics[width=4.8cm]{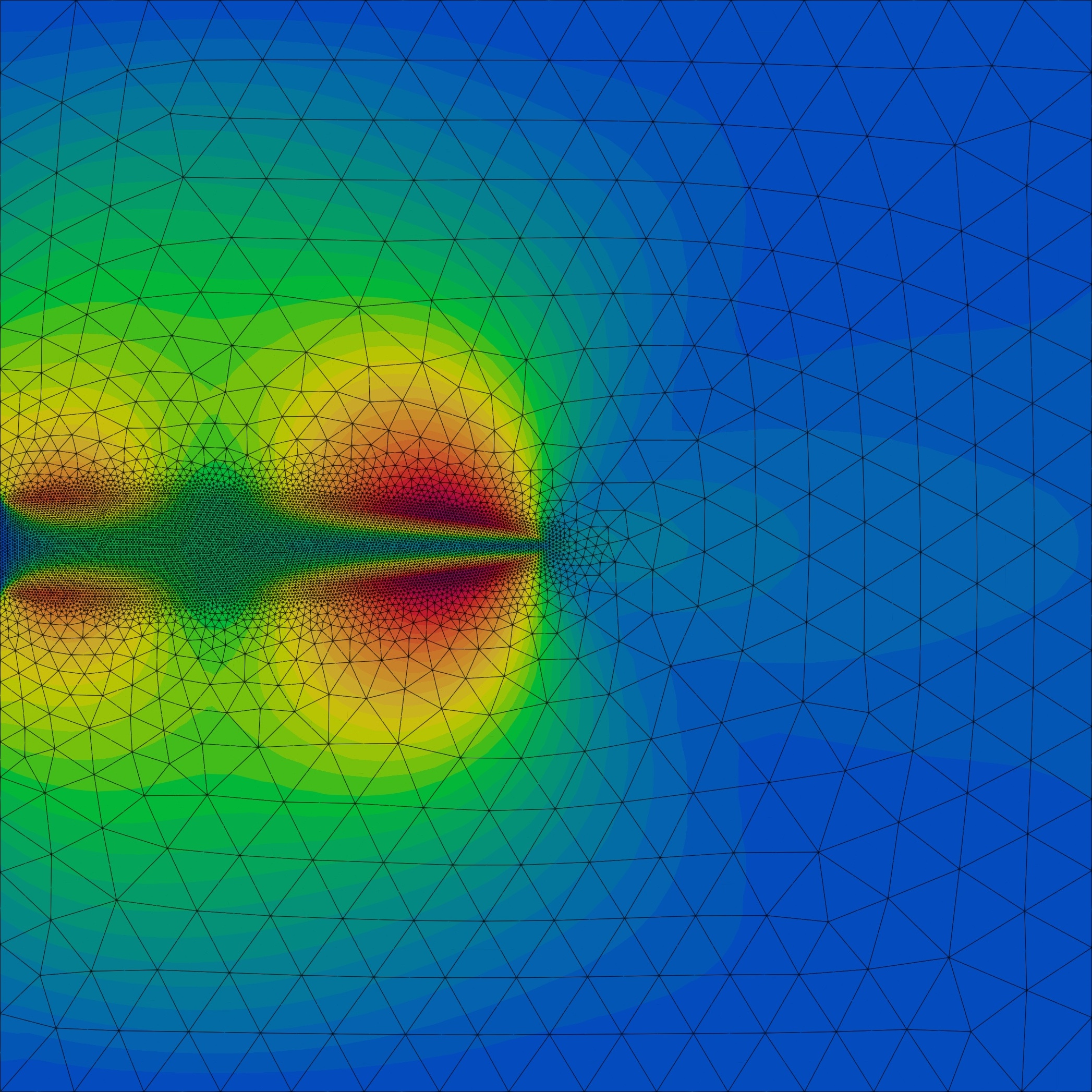}
    \vskip.5\baselineskip
    \includegraphics[width=4.8cm]{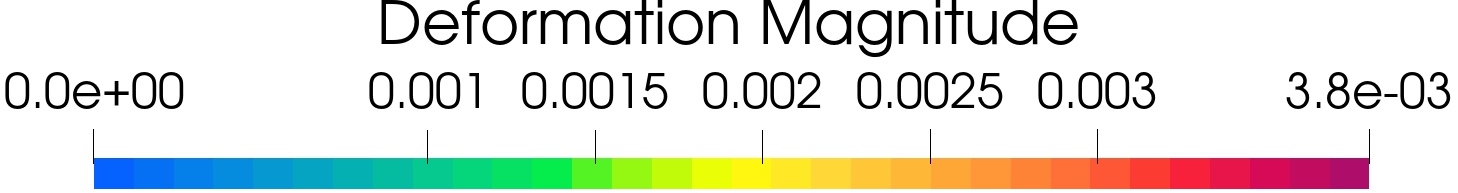}
  \end{minipage}
  \begin{minipage}[t]{.32\textwidth}
    \centering
    \includegraphics[width=4.8cm]{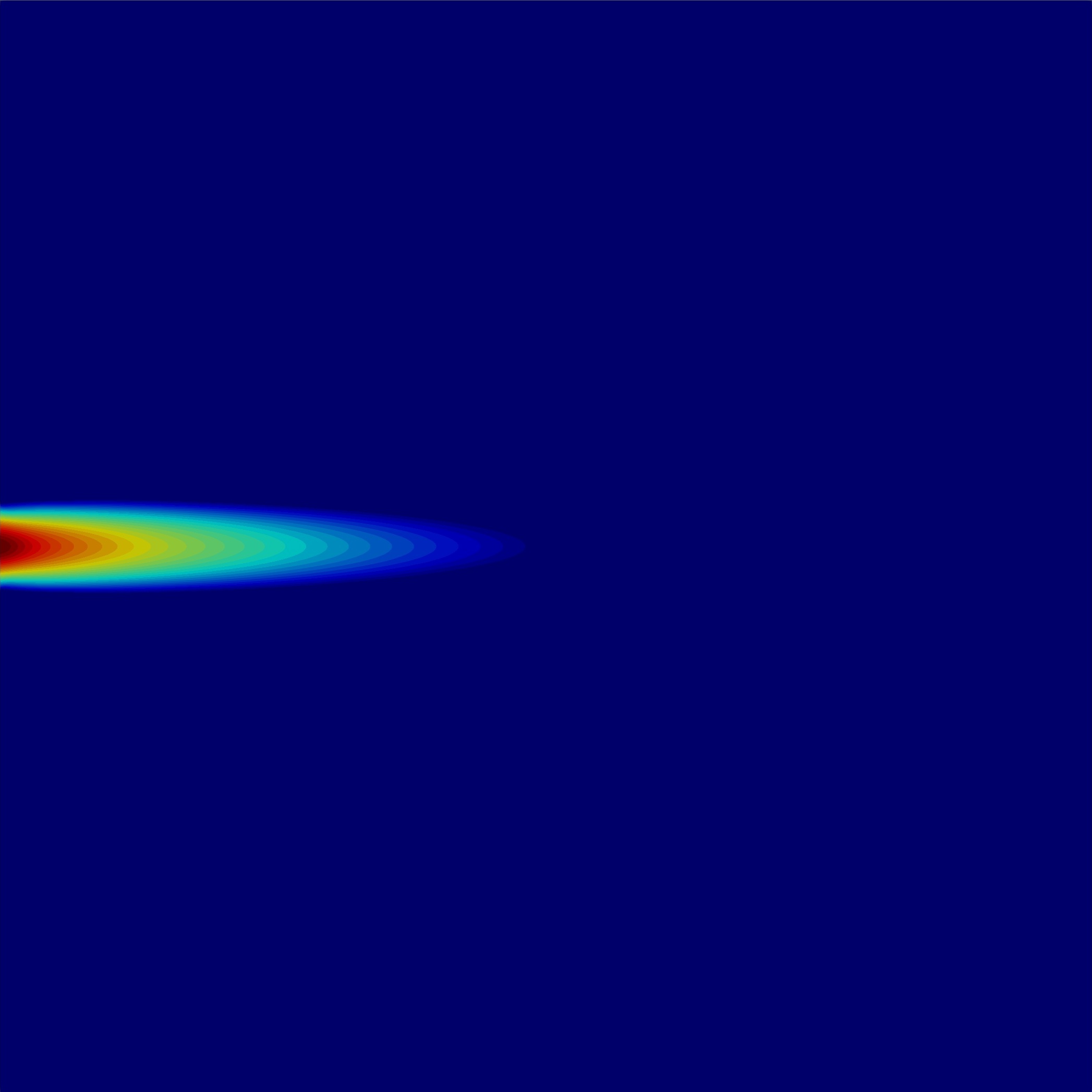}
    \vskip.5\baselineskip
    \includegraphics[width=4.8cm]{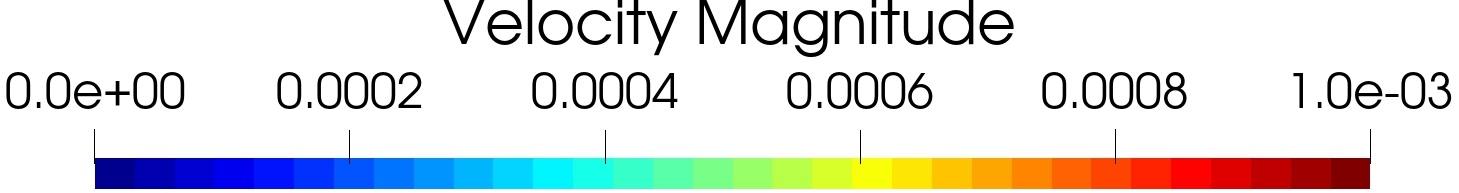}
  \end{minipage}
  \begin{minipage}[t]{.32\textwidth}
    \centering
    \includegraphics[width=4.8cm]{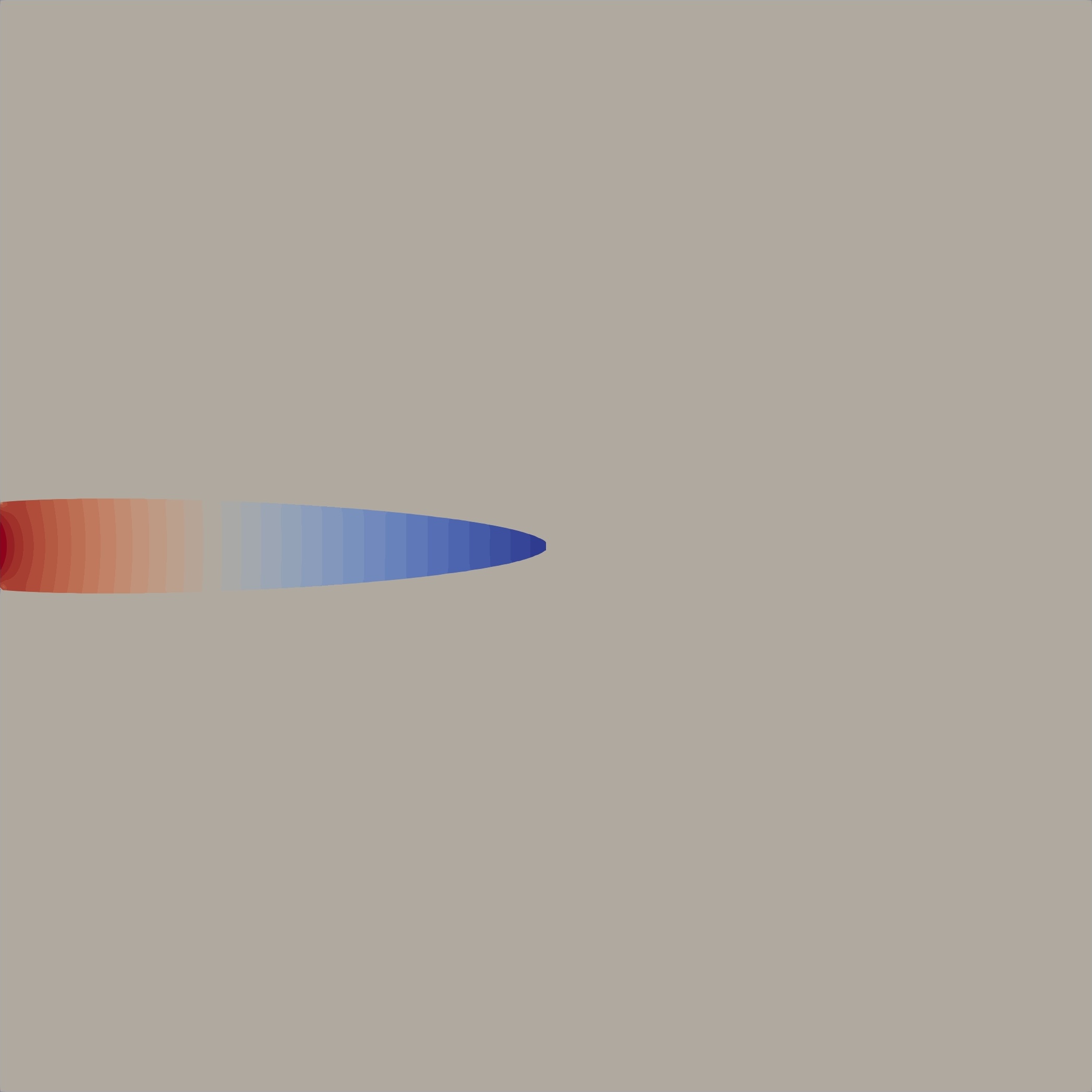}
    \vskip.5\baselineskip
    \includegraphics[width=4.8cm]{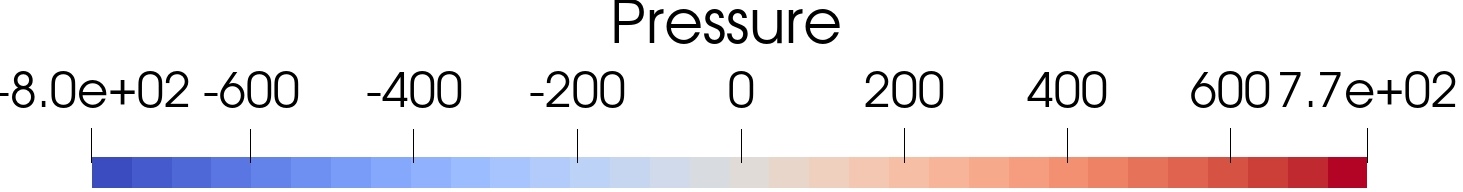}
  \end{minipage}
  \caption{Example 3: FSI-deformation magnitude, velocity magnitude and FSI-pressure after three levels of mesh refinement.}
  \label{fig.results.ex3:fsi}
\end{figure}

\begin{figure}
  \centering
  \includegraphics{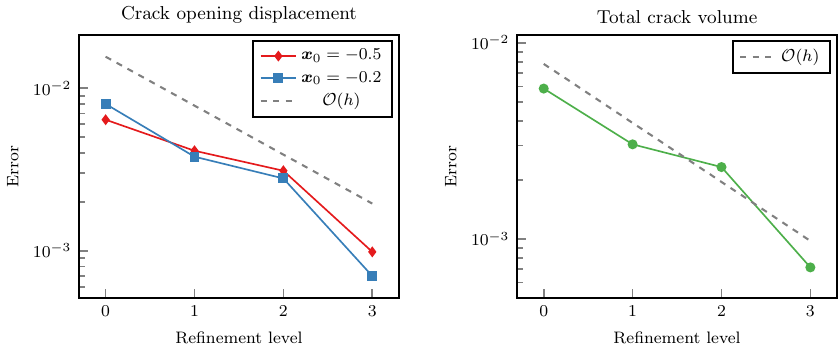}
  \caption{Example 3: Convergence of $\COD$ and $\TCV$ on different spatial refinement levels.}
  \label{fig.results.ex3:convergence}
\end{figure}

We consider $h_0 = 0.01$ and five levels of mesh refinement. The phase-field deformation, together with the computational mesh, and the phase-field solution on mesh level three, can be seen in \Cref{fig.results.ex3:phase-field}. In particular, we see the expected discontinuity of the vertical component of the deformation along the crack. In \Cref{fig.results.ex3:fsi}, we observe the fluid-structure interaction deformation, velocity and pressure on the same mesh level. This particularly illustrates the reconstructed domain and the flexibility regarding the problem posed inside the open crack by the domain reconstruction. Finally, In \Cref{fig.results.ex3:convergence}, we again consider the numerical convergence of the $\COD$ and $\TCV$ towards the values resulting from the finest mesh level. We again observe linear convergence, with a slight deviation on mesh level two.

\subsection{Example 4: A propagating fracture}
We extend the previous examples by considering an incrementally-increasing,
pressure resulting in a growing crack.

\subsubsection{Configuration}
The domain of interest is again $\Omega=(-2, 2)^2$. The mechanical parameters 
are as in Example~1, that is, $E_s = 10^5$ and $\nu_s = 0.35$. We choose the 
model parameters in the phase-field problem as $G_c=500$, $\kappa=10^{-10}$, 
$\gamma = 100 h^{-2}$, while we take the phase-field regularisation parameter 
fixed as $\varepsilon=0.04$. This is because we cannot expect mesh convergence 
in this setting with varying $\varepsilon$ \Cite{HWW15,MiWheWi14}.

The time-dependent background pressure driving the growing crack is
\begin{equation*}
  p(t) = 4500 + 10t \qquad\text{for }t\in[0, 0.5].
\end{equation*}
We take the constant time step $\Delta t = 0.1$. In each time steps, we again
use five iterations so solve the phase-field problem.

For the fluid-structure interaction problem, we take the right-hand side 
$\fbref=(5, 0)^T$. In a pure fluid flow problem, this would be a no-flow 
problem. But even in the FSI context, this results in a large pressure. We 
chose this as only the pressure couples back to the phase-field  problem, 
and we need a sufficiently large pressure to have an impact on the dynamics 
of the phase-field problem. With our given right-hand side forcing 
term, the pressure difference between the left an right tip is about 
$2\cdot 10^3$.

\subsubsection{Results}
\begin{figure}
  \centering
  \includegraphics[width=5cm, trim=92.4cm 110.8cm 92.4cm 110.8cm, clip]{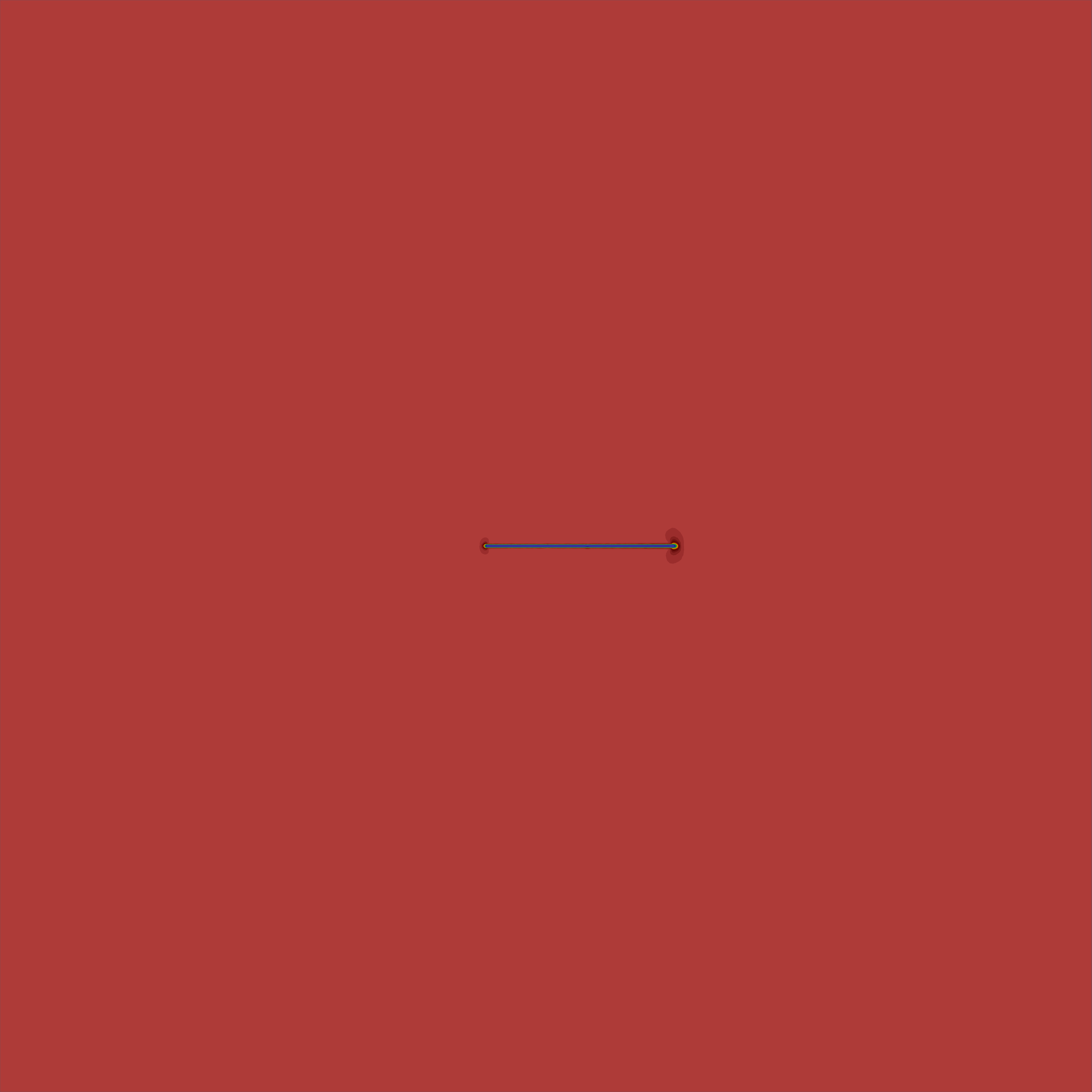}
  \includegraphics[width=5cm, trim=92.4cm 110.8cm 92.4cm 110.8cm, clip]{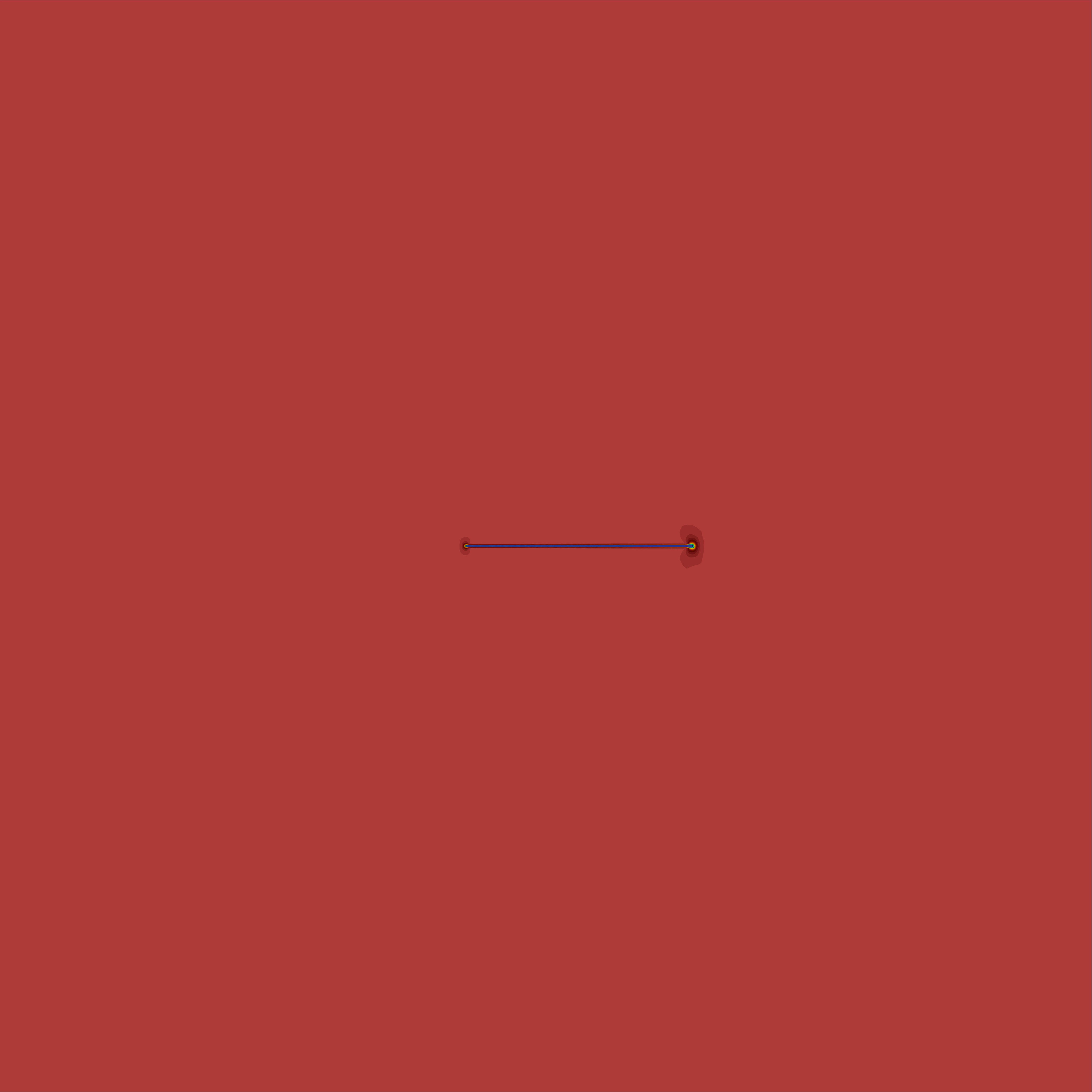}
  \includegraphics[width=5cm, trim=92.4cm 110.8cm 92.4cm 110.8cm, clip]{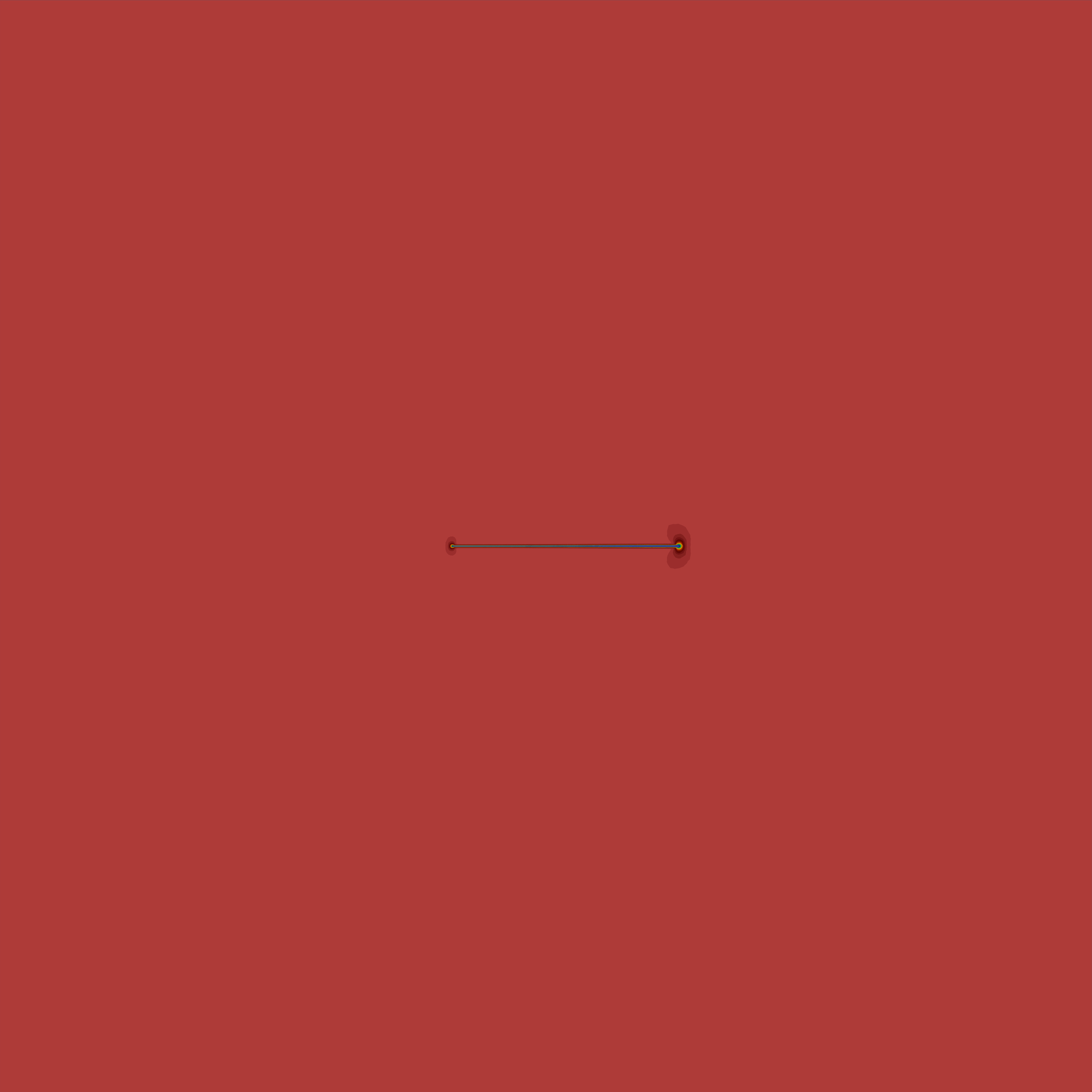}
  \caption{Example 4: Phase-field solution at $t=0.5$ over three meshes in the region $(-0.75,0.75)\times(-0.5,0.5)$.}
  \label{fig.example4.phase-field}
\end{figure}

\begin{figure}
  \centering
  \includegraphics{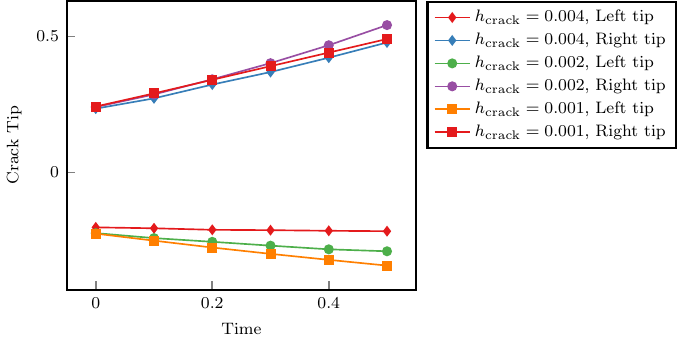}
\caption{Example 4: Crack tips over time.}
  \label{fig.example4.tips}
\end{figure}

The phase-field solution on meshes with $h_\text{crack}=0.004, 0.002, 0.001$ 
is shown in \Cref{fig.example4.phase-field}. The temporal behaviour of the 
crack tip is plotted in \Cref{fig.example4.tips}. Furthermore, at $t=1$, the
left mesh tip is at $x=-0.21625492$, $-0.29008647$, $-0.34273292$ and the right 
tip at $x=0.47777961$, $0.54183778$, $0.49064490$ on the three meshes 
respectively. As expected, the crack grows faster towards the right, where 
the total pressure is larger than on the left due to the coupled FSI pressure. 
While the crack appears to grow only marginally on the coarsest mesh, the 
results are consistent on the two finer meshes.

\section{Conclusions}
\label{sec_conclusions}
In this work, we designed and tested a fully coupled algorithm for iterating 
between phase-field fracture and fluid-structure interaction. More 
specifically, we first designed a pressurized phase-field approach 
including the interface conditions stated explicitly on the boundary of
the crack domain to compute the fracture path or fracture propagation. 
This is motivated by the fact that the fluid-structure interaction pressure,
used to couple the FSI and phase-field models, only exists inside the
fluid-filled crack and not the entire domain.

Next, we used the fracture displacements resulting from the phase-field model 
to reconstruct the mesh to distinguish between the fluid-filled fracture domain 
and the surrounding medium by a sharp interface. With the reconstructed 
mesh, we used the well-known interface-tracking approach (ALE) to solve the 
fluid-structure interaction problem. This allowed us to solve, for example, 
Stokes flow in the fracture domain, driven by both volumetric forcing or 
prescribed boundary conditions. The resulting pressure is applied to 
the pressurized phase-field fracture sub-problem, which allowed us to couple 
the two sub-problems iteratively.

We substantiated our final algorithm with numerical experiments in which 
we undertook both quantitative studies, in terms of computational convergence 
studies, and observed qualitative behaviour of the solution. The
current findings are promising and give hope to extend this in further 
studies to fully time-dependent problems. However, for the latter, the 
entire model must be updated, including time-dependent terms, which 
are planned as future work.

\section*{Data Availability Statement}
The source code and the data generated with it is publicly available on 
github under \url{https://github.com/hvonwah/repro-coupled-phase-field-fsi}
and archived on zenodo~\cite{vWW23_zenodo}.

\section*{Acknowledgments}
HvW acknowledges support through the Austrian Science Fund (FWF) project F65.

TW thanks Ivan Yotov (University of Pittsburgh) for a fruitful discussion at the 
`Hot Topics: Recent Progress in Deterministic and Stochastic Fluid-Structure Interaction' Workshop, Dec 2023, in Berkeley, US, during the revision of this paper.

\printbibliography

\end{document}